\documentclass[review,onefignum,onetabnum]{siamart171218}

\usepackage{lipsum}
\usepackage{amsfonts}
\usepackage{graphicx}
\usepackage{epstopdf}
\usepackage{algorithmic}
\usepackage{color}
\usepackage{mathrsfs}
\usepackage{amssymb}
\usepackage{amsmath}
\usepackage{boxedminipage}
\usepackage{stmaryrd}
\usepackage{multirow}
\usepackage{booktabs}
\usepackage{accents}
\ifpdf
  \DeclareGraphicsExtensions{.eps,.pdf,.png,.jpg}
\else
  \DeclareGraphicsExtensions{.eps}
\fi

\newcommand{\RR}{\mathbf R}
\newtheorem{assumption}{Assumption}
\newsiamremark{remark}{Remark}
\newsiamremark{hypothesis}{Hypothesis}
\crefname{hypothesis}{Hypothesis}{Hypotheses}
\newsiamthm{claim}{Claim}

\headers{Stochastic Primal-Dual Coordinate Method for NCCP}{Daoli Zhu and Lei Zhao}

\title{Stochastic Primal-Dual Coordinate Method for Nonlinear Convex Cone Programs\thanks{Submitted to the editors DATE.
\funding{This work was funded NSFC: 71471112.}}}

\author{Daoli Zhu\thanks{Antai College of Economics and Management and Sino-US Global Logistics Institute, Shanghai
Jiao Tong University, 200030 Shanghai, China
  (\email{dlzhu@sjtu.edu.cn})}
\and Lei Zhao\thanks{Antai College of Economics and Management and Sino-US Global Logistics Institute, Shanghai
Jiao Tong University, 200030 Shanghai, China
  (\email{l.zhao@sjtu.edu.cn})}  
}

\usepackage{amsopn}

\ifpdf
\hypersetup{
  pdftitle={Stochastic Primal-Dual Coordinate Method for Nonlinear Convex Cone Programs},
  pdfauthor={Daoli Zhu and Lei Zhao}
}
\fi

\begin{document}

\maketitle

\begin{abstract}
Block coordinate descent (BCD) methods and their variants have been widely used in coping with large-scale nonconstrained optimization problems in many fields such as imaging processing, machine learning, compress sensing and so on. For problem with coupling constraints, Nonlinear convex cone programs (NCCP) are important problems with many practical applications, but these problems are hard to solve by using existing block coordinate type methods. This paper introduces a stochastic primal-dual coordinate (SPDC) method for solving large-scale NCCP. In this method, we randomly choose a block of variables based on the uniform distribution. The linearization and Bregman-like function (core function) to that randomly selected block allow us to get simple parallel primal-dual decomposition for NCCP. The sequence generated by our algorithm is proved almost surely converge to an optimal solution of primal problem. Two types of convergence rate with different probability (almost surely and expected) are also obtained. The probability complexity bound is also derived in this paper.
\end{abstract}



\section{Introduction}\label{Intro}
Recently, block coordinate descent (BCD) methods and their variants have been widely used in coping with large-scale optimization problems in many fields such as imaging processing, machine learning, compress sensing and so on (See~\cite{LiOsher, QinScheinberg, WuLange, YunToh}). In each iteration of these methods, one block of coordinates is selected to reduce the objective value, simultaneously the other blocks are keeping fixed. The main advantage of BCD is to reduce the complexity and memory requirements per iteration. These benefits are increasingly important for very-large scale problem.\\
\indent For problem without constraints, there are two variations of BCD discussed the most by researchers. The first variation is on block choosing strategy. One common approach for block choosing is cyclic strategy. Tseng~\cite{Tseng} proved the convergence of a BCD of cyclic strategy. Luo and Tseng~\cite{LuoZeng} and Wang and Lin~\cite{WangLin} proved local and global linear convergence under specific assumptions respectively. The other approach is randomized strategy. Nesterov~\cite{Nesterov} studied the convergence rate of randomized BCD for convex smooth optimization. Richt\'arik and Tak\'a\v c~\cite{RT} and Lu and Xiao~\cite{LuXiao} extend Nesterov's technique to composite optimization. The point read to evaluate the gradient in each iteration is the second variation of BCD. If the read point have different "ages", this type of BCD called asynchronous BCD, otherwise it is called synchronous BCD. All the variants of BCD reviewed above are synchronous BCD. Liu and Wright~\cite{LiuWright1} and Liu et. al.~\cite{LiuWright2} established the convergence rate of asynchronous BCD for composite optimization and convex smooth optimization without constraints respectively.\\
\indent For problem with constraints, there are only a few works. Gao et. al.~\cite{GaoXuZhang} proposed a coordinate type method for problem with linear coupling constraints. Necoara and Patrascu~\cite{Necoara} proposed a random coordinate descent algorithm for optimization problem with one linear constraint. To the best of our knowledge, there is no primal-dual coordinate convergence and convergence rate results for optimization problem with nonlinear inequality constraints.\\
\indent In the big data era, there is one class of optimization problem with nonlinear inequality constraints which has received much attention in the compressed sensing and machine learning fields, such as least absolute shrinkage and selection operator (LASSO) problem, basis pursuit denoising (BPDN) problem and support vector machine (SVM) problem. Actually, the LASSO problem is a least-squares problem with a $1$-norm constraint (Tibshirani~\cite{LASSO}) and BPDN is a sparse optimization problem with a least-squares constraint (Chen et. al.~\cite{ChenDonohoSaunders}). In the field of SVM, Oneto et. al.~\cite{SVM} states that the Ivanov regularization scheme allows more effective control of the learning hypothesis space, while the Morozov regularization scheme is useful when an effective estimate of the noise in the data is available. All these problems are problem with one smooth or nonsmooth inequality constraint.\\
\indent In this paper, we consider the following nonlinear convex cone programs (NCCP):
\begin{equation}\label{Prob:general-function}
\begin{array}{lll}
\mbox{(P):}  &\min       & G(u)+J(u)      \\
             &\rm {s.t}  & \Theta(u)\in -\mathbf{C} \\
             &           & u\in \mathbf{U}.
\end{array}
\end{equation}
\noindent where $G$ is a convex smooth function on closed convex subset $\mathbf{U}\subset \RR^{n}$, $J$ is a convex but possibly nonsmooth function on $\mathbf{U}\subset \RR^{n}$, and $\Theta(u)$ is a possibly nonsmooth $\mathbf{C}$-convex mapping from $\RR^{n}$ to $\RR^{m}$. $\mathbf{C}$ is a nonempty closed convex cone in $\RR^{m}$ with vertex at the origin, that is $\alpha\mathbf{C}+\beta\mathbf{C}\subset \mathbf{C}$, for $\alpha,\beta\geq 0$. It is obvious that when $\mathring{\mathbf{C}}$ (the interior of $\mathbf{C}$) is nonempty, the constraint $\Theta(u)\in -\mathbf{C}$ corresponds to an inequality constraint. The case $\mathbf{C}=\{0\}$ corresponds to an equality constraint. $\mathbf{C}^{*}$ denotes the conjugate cone. Assume that both $J(u)=\sum\limits_{i=1}\limits^{N}J_i(u_i)$ and $\Theta(u)=\sum\limits_{i=1}\limits^{N}\Theta_{i}(u_{i})$ are additive respect to following space decomposition:
\begin{equation}\label{spacedecomposition}
\mathbf{U}=\mathbf{U}_1\times\mathbf{U}_2\cdots\times\mathbf{U}_N, u_i\in\mathbf{U}_i\subset \RR^{n_i}\qquad\mbox{and}\qquad \sum\limits_{i=1}\limits^{N}n_i=n.
\end{equation}
\indent In order to cope with separable or nonseparable, smooth or nonsmooth NCCP, Cohen and Zhu~\cite{CohenZ} proposed a fairly general first-order primal-dual decomposition method, Auxiliary Problem Principle of augmented Lagrangian methods (APP-AL), based on linearization of the augmented Lagrangian. Recently, Zhao and Zhu~\cite{ZhaoZhu2017} extend APP to a more flexible framework, Variant Auxiliary Problem Principle methods (VAPP), with an O(1/t) convergence rate on average for primal suboptimality, feasibility, and dual suboptimality. Moreover, they propose a backtracking scheme to treat the case where the Lipschitz constants are not known or computable.\\
\indent In this paper, we propose a stochastic primal-dual coordinate (SPDC) method based on auxiliary problem principle (Cohen and Zhu~\cite{CohenZ}) for NCCP. In this method, we randomly choose a block of variables based on the uniform distribution. The linearization and Bregman-like function (core function) to that randomly selected block allow us to get simple parallel primal-dual decomposition for problem (P). The sequence generated by our algorithm is proved almost surely converge to an optimal solution of problem (P). The almost surely $O(1/t^{\frac{1-\alpha}{2}})$ convergence rate and expected $O(1/t^{1-\alpha})$ result with $1/2<\alpha<1$ are also obtained for problem (P) under the convexity assumptions. The probability complexity bound is also derived in this paper.\\
\indent The rest of this paper is organized as follows. Section~\ref{Pre} is devoted to technical preliminaries. The updating scheme of SPDC for (P) is presented in Section~\ref{SPDC}. In Section~\ref{Convergence} we establish the convergence. Almost surely $O(1/t^{\frac{1-\alpha}{2}})$ sub-linear convergence rate is proposed in Section~\ref{ARate}. In Section~\ref{ERate} expected $O(1/t^{1-\alpha})$ sub-linear convergence rate and the probability complexity bound is established. In Section~\ref{sec:experiments}, we compute an example to show the efficiency of our algorithm.
\section{Preliminaries}\label{Pre}
In this section, we first provide some preliminaries which are useful for our further discussions and then summarize some notations and assumptions to be used. We denote $\langle \cdot \rangle $ and $\| \cdot \|$ as the inner product and Euclidean norm of vector, respectively.
\subsection{Notations and assumptions}
Throughout this paper, we make the following standard assumptions for Problem (P):
\begin{assumption}\label{assump1}
{\rm
\noindent \begin{itemize}
\item[(i)] $J$ is a convex, l.s.c function such that $\mathbf{dom}J\cap \mathbf{U}\neq\emptyset$, $J$ is not necessary differentiable. $J$ is subgradientiable and has linear bounded subgradients in $\mathbf{U}$, that is
    \begin{eqnarray*}
    \exists c_1>0, c_2>0,\;\forall u\in\mathbf{U},\;\forall r\in\partial J(u),\;\|r\|\leq c_1\|u\|+c_2.
    \end{eqnarray*}
\item[(ii)] $G$ is a convex and differentiable with its derivative Lipschitz of constant $B_{G}$.
\item[(iii)] $G+J$ is coercive on $\mathbf{U}$, if $\mathbf{U}$ is not bounded, that is
\begin{eqnarray*}
\forall \{u^{k}|k\in \mathbb{N}\}\subset \mathbf{U}, \lim_{k\rightarrow+\infty}\|u^{k}\|=+\infty\Rightarrow\lim_{k\rightarrow+\infty}(G+J)(u^{k})=+\infty.
\end{eqnarray*}
\item[(iv)] $\Theta$ is $\mathbf{C}$-convex, where
\begin{equation}\label{Theta_C_convex}
\forall u,v\in \mathbf{U},\;\forall \alpha\in [0,1],\;\Theta(\alpha u+(1-\alpha)v)-\alpha\Theta(u)-(1-\alpha)\Theta(v)\in - \mathbf{C}.
\end{equation}
Moreover, $\Theta(u)$ is Lipschitz with constant $\tau$ on an open subset $\mathcal{O}$ containing $\mathbf{U}$, where
\begin{equation}\label{Theta_Lipschitz}
\forall u,v\in \mathcal{O},\;\|\Theta(u)-\Theta(v)\|\leq\tau\|u-v\|.
\end{equation}
\item[(v)] Constraint Qualification Condition. When $\mathring{\mathbf{C}}\neq\emptyset$, we assume that
\begin{equation}\label{Theta_CQC_neq}
\Theta(\mathbf{U})\cap(-\mathring{\mathbf{C}})\neq\emptyset.
\end{equation}
\end{itemize}
}
\end{assumption}
Furthermore, the following lemma gives some properties of $G$ which satisfies Assumption~\ref{assump1}.
\begin{lemma}\label{lemma:fphi}
Suppose function $G$ satisfy (ii) of Assumption~\ref{assump1} on $\mathbf{U}$, then we have that:
for all $u,v\in \mathbf{U}$, $\langle\nabla G(v),v-u\rangle\leq G(v)-G(u)+\frac{B_G}{2}\|u-v\|^2$.
\end{lemma}
\indent The results of this lemma are classical, the proof is omitted. (see Zhu and Marcotte~\cite{Zhu96})
\subsection{Augmented Lagrangian and first-order primal-dual decomposition algorithm}
In this subsection, the Lagrangian of (P) is defined as:
\begin{equation}\label{func:L}
L(u,p)=(G+J)(u)+\langle p,\Theta(u)\rangle,
\end{equation}
and a saddle point $(u^*,p^*)\in \mathbf{U}\times \mathbf{C}^*$ is such that
\begin{equation}
\forall u\in \mathbf{U},\; \forall p\in \mathbf{C}^{*}:\; L(u^{*},p)\leq L(u^{*},p^{*})\leq L(u,p^{*}). \label{saddle point:L}
\end{equation}
\indent Under Assumption~\ref{assump1}, there exist saddle points of $L$ on $\mathbf{U}\times \mathbf{C}^*$. The dual function $\psi$ is defined as
\begin{equation*}
\psi(p)=\begin{cases}\displaystyle\min_{u\in \mathbf{U}}L(u,p)  &\forall p\in\mathbf{C}^*\\
                     \displaystyle-\infty  &\mbox{otherwise.}\end{cases}
\end{equation*}
\indent The function $\psi$ is concave and sub-differentiable. Using dual function $\psi(p)$, we consider the primal-dual pair of nonlinear convex cone optimization:
\begin{equation*}
\begin{array}{lllllll}
  \mbox{(P):} &\min      & (G+J)(u)                 &\qquad\qquad\qquad\qquad&\mbox{(D):}    &\max      & \psi(p) \\
              &\rm {s.t} & \Theta(u)\in -\mathbf{C} &\qquad\qquad\qquad\qquad&               &\rm {s.t} & p\in\mathbf{C}^*\\
              &          & u\in\mathbf{U}           &\qquad\qquad\qquad\qquad&               &          &
\end{array}
\end{equation*}
\indent The following theorem characterizes a saddle point optimality condition for the primal and dual problem.
\begin{theorem}\label{theo:Lag}
A solution $(u^*,p^*)$ with $u^*\in \mathbf{U}$ and $p^*\in \mathbf{C}^*$ is a saddle point for the Lagrangian function $L(u,p)$ if and only if
\begin{itemize}
\item[{\rm(i)}] $L(u^*,p^*)=\min\limits_{u\in \mathbf{U}}L(u,p^*)$\\
           or the following variational inequality holds:$$\langle\nabla G(u^*),u-u^*\rangle+J(u)-J(u^*)+\langle p^*,\Theta(u)-\Theta(u^*)\rangle\geq 0, \forall u\in\mathbf{U};$$
\item[{\rm(ii)}] $\Theta(u^*)\in -\mathbf{C}$;
\item[{\rm(iii)}] $\langle p^*,\Theta(u^*)\rangle=0$.
\end{itemize}
Moreover, $(u^*,p^*)$ is a saddle point if and only if $u^*$ and $p^*$ are, respectively, optimal solutions to the primal and dual problems (P) and (D) with no duality gap, that is, with $(G+J)(u^*)=\psi(p^*)$.
\end{theorem}
\indent Now we take a trick by introducing slack variables which help problem (P) come back to problem with equality constraints. Namely, the problem (P) is converted into the equivalent problem with equality constraints as follows
\begin{equation*}
\begin{array}{lllllll}
  \mbox{(P$_1$):} &\min\limits_{\xi\in-\mathbf{C}}      & (G+J)(u) \\
              &\rm {s.t} & \Theta(u)-\xi=0 \\
              &          & u\in\mathbf{U}  \\
\end{array}
\end{equation*}
\indent The augmented Lagrangian for this problem is
\begin{equation}\label{L_xi}
\overline{L}_\gamma(u,\xi,p)=(G+J)(u)+\langle p,\Theta(u)-\xi\rangle+\frac{\gamma}{2}\|\Theta(u)-\xi\|^2
\end{equation}
\indent The augmented Lagrangian associated with problem (P) is defined as
\begin{equation}\label{L_gamma}
L_\gamma(u,p)\triangleq\min_{\xi\in-\mathbf{C}}\overline{L}_\gamma(u,\xi,p)=(G+J)(u)+\varphi(\Theta(u),p),
\end{equation}
where $\varphi(\Theta(u),p)=[\|\Pi\big{(}p+\gamma\Theta(u)\big{)}\|^2-\|p\|^2]/2\gamma$ and $\Pi$ is a projection on to $\mathbf{C}^*$.\\
\indent The augmented Lagrangian dual function is as following:
\begin{equation}
\forall p\in \RR^m, \psi_\gamma(p)=\min_{u\in\mathbf{U}}L_\gamma(u,p)=\min_{u\in\mathbf{U}}(G+J)(u)+\varphi(\Theta(u),p).
\end{equation}
\indent Using $\psi_\gamma(p)$, we obtain new primal-dual pair of nonlinear convex cone optimization
\begin{equation*}
\begin{array}{lllllll}
  \mbox{(P):} &\min      & (G+J)(u)                 &\qquad\qquad\qquad\qquad&\mbox{(D$_\gamma$):}    &\max      & \psi_\gamma(p) \\
              &\rm {s.t} & \Theta(u)\in -\mathbf{C} &\qquad\qquad\qquad\qquad&                      &\rm {s.t} & p\in\mathbf{R}^m\\
              &          & u\in\mathbf{U}           &\qquad\qquad\qquad\qquad&                      &          &
\end{array}
\end{equation*}
\indent The following theorem shows that function $\varphi(\theta,p)$, dual function $\psi_\gamma(p)$ and augmented Lagrangian $L_\gamma(u,p)$ have some useful properties.
\begin{theorem}\label{theo:Lag_2}
Suppose Assumption~\ref{assump1} holds for problem (P). Then we have
\begin{itemize}
\item[{\rm(i)}] The function $\varphi(\theta,p)$ is convex in $\theta$ and concave in $p$.
\item[{\rm(ii)}] $\varphi$ is differentiable in $\theta$ and $p$ and one has
\begin{eqnarray*}
&&\nabla_{\theta}\varphi(\theta, p)=\Pi(p+\gamma\theta),\\
&&\nabla_{p}\varphi(\theta, p)=[\Pi(p+\gamma\theta)-p]/\gamma,\\
&&\varphi(\theta, p)=[\|\Pi(p+\gamma\theta)\|^2-\|p\|^2]/2\gamma.
\end{eqnarray*}
\item[{\rm(iii)}] $\psi_\gamma(p)$ is concave and differentiable in $p$, and $\nabla\psi_\gamma(p)=[\Pi(p+\gamma\Theta(\hat{u}(p)))-p]/\gamma$, where $\hat{u}(p)\in\hat{\mathbf{U}}(p)=\{u\in\mathbf{U}|u=\arg\min\limits_{u\in\mathbf{U}}L_\gamma(u,p)\}$.
\item[{\rm(iv)}] $L$ and $L_\gamma$ have the same sets of saddle points $\mathbf{U}^*\times\mathbf{P}^*$ respectively on $\mathbf{U}\times\mathbf{C}^*$ and $\mathbf{U}\times\mathbf{R}^m$.
\item[{\rm(v)}] $L_\gamma$ is stable in $u$, that is $\forall p^*\in \mathbf{P}^*, \hat{\mathbf{U}}(p^*)=\mathbf{U}^*$.
\end{itemize}
\end{theorem}
\indent For the general NCCP, the augmented Lagrangian method is an approach which can overcome the instability and nondifferentiability of the dual function of the Lagrangian. Furthermore, the augmented Lagrangian of a constrained convex program has the same solution set as the original constrained convex program. The augmented Lagrangian approach for equality-constrained optimization problems was introduced in Hestenes~\cite{Hestenes1969} and Powell~\cite{Powell1969}, then extended to inequality-constrained problems by Buys~\cite{Buys1972}.\\
\indent Although the augmented Lagrangian approach (Uzawa algorithm) has several advantages, it does not preserve separability, even when the initial problem is separable. One way to decompose the augmented Lagrangian is ADMM (Fortin and Glowinski~\cite{ADMM1983}). ADMM can only handle convex problems with linear constraints and is not easily parallelizable. Another way to overcome this difficulty is the Auxiliary Problem Principle of augmented Lagrangian methods (APP-AL) (Cohen and Zhu~\cite{CohenZ}), which is a fairly general first-order primal-dual decomposition method based on linearization of the augmented Lagrangian in separable or nonseparable, smooth or nonsmooth NCCP. Cohen and Zhu (1984)~\cite{CohenZ} propose first-order primal-dual augmented Lagrangian methods for NCCP as Algorithm~\ref{alg:APP-AL}.\\
\begin{algorithm}
\caption{Auxiliary Problem Principle for Augmented Lagrangian (APP-AL)}\label{alg:APP-AL}
\begin{algorithmic}
\STATE{Initialize $u^0 \in \mathbf{U}$ and $p^0\in \mathbf{C^*}$}
\FOR{$k = 0,1,\cdots $,}
\STATE{Update $u^{k+1}=\arg\min\limits_{u\in \mathbf{U}}\langle\nabla G(u^{k}), u \rangle + J(u)+ \langle\Pi\big{(}p^k+\gamma\Theta(u^k)\big{)}, \Theta(u)\rangle+\frac{1}{\epsilon}D(u,u^k)$}
\STATE{Update $p^{k+1}=\Pi\big{(}p^k+\gamma\Theta(u^{k+1})\big{)}$}
\ENDFOR
\end{algorithmic}
\end{algorithm}
\indent Cohen and Zhu (1984) shows the sequence $\{(u^k,p^k)\}$ generated by APP-AL convergence to $(u^*,p^*)$ saddle point of $L$ over $\mathbf{U}\times\mathbf{C}^*$. If we choose the following additive $K(u)$ respect to the space decomposition~\ref{spacedecomposition}, $K(u)=\sum\limits_{i=1}\limits^{N}K_i(u_i)$,  then the primal problem split to $N$ small optimization subproblem only for $u_i$.~\cite{CohenZ} In the era of big data, there has been a surge of interest in redesign of APP-AL suitable for solving the huge optimization with available computing performance.
\subsection{The properties of projection on convex cone}
In this subsection, we introduce some properties of projection on convex sets (resp. convex cone) as preparations. These properties are used in the following sections.\\
\indent Let $\mathcal{S}$ be a nonempty closed convex set of $\RR^m$. For $x\in\RR^m$, we propose the projection $\Pi_{\mathcal{S}}(x)$ as a projection on $\mathcal{S}$. Then $\Pi_{\mathcal{S}}(x)$ is characterized by the following two conditions~\cite{Projectiononconvexsets}:
\begin{eqnarray}
&(i)&\qquad\langle y-\Pi_{\mathcal{S}}(x), x-\Pi_{\mathcal{S}}(x)\rangle\leq0, \forall y\in\mathcal{S};\label{eq:Projecproperty3}\\
&(ii)&\qquad\|\Pi_{\mathcal{S}}(x)-\Pi_{\mathcal{S}}(y)\|\leq\|x-y\|.\label{eq:Projecproperty4}
\end{eqnarray}
\indent Furthermore, the following proposition gives another property of projection operator which is used for convergence and convergence rate analysis.
\begin{proposition}\label{proposition}
\noindent For any $(x,y,z)\in\RR^{m\times m\times m}$, the projection operator $\Pi_{\mathcal{S}}$ satisfies
\begin{equation}\label{eq:Projecproperty7}
2\langle\Pi_{\mathcal{S}}(z+x)-\Pi_{\mathcal{S}}(z+y),x\rangle\leq\|x-y\|^2+\|\Pi_{\mathcal{S}}(z+x)-z\|^2-\|\Pi_{\mathcal{S}}(z+y)-z\|^2.
\end{equation}
\end{proposition}
\begin{proof} See~\cite{ZhaoZhu2017}.
\end{proof}
\indent Next, we consider the properties for projection on convex cone. Let $\mathbf{C}$ be a nonempty closed convex cone in $\mathbf{R}^m$ with vertex at the origin. $\mathbf{C}^{*}$ denotes the conjugate cone. Let $\Pi$ denote the projection on $\mathbf{C}^*$ and $\Pi_{-\mathbf{C}}$ denote the projection on $-\mathbf{C}$. The projection is characterized by the following conditions. (see Wierzbicki~\cite{Projection}):
\begin{eqnarray}
&(iii)  &\qquad y=\Pi(y)+\Pi_{-\mathbf{C}}(y), y\in\mathbf{R}^m;\label{eq:Projecproperty5}\\
&(iv) &\qquad \langle\Pi(y), \Pi_{-\mathbf{C}}(y)\rangle=0, y\in\mathbf{R}^m.\label{eq:Projecproperty6}
\end{eqnarray}
\subsection{The boundness of multiplier sets and application}
Noted that under Assumption~\ref{assump1}, the dual optimal solution set of Lagrangian for (P) is bounded~\cite{CohenZ}. Here we denote the bound is $\mu_0$ i.e., $\|p^{*}\|\leq\mu_0$. Moreover, we set $\mu=\mu_0+1$. The estimate of dual optimal bound $\mu$ is required in constructing new coordinate descent method. In this section, we will provide the estimate of dual optimal bound for problem (P) with special convex cone $\mathbf{C}=\RR_+^m$ or $\mathbf{C}=\mathcal{K}_{\nu}^m$.\\
\indent If $\mathbf{C}=\RR_+^m$ and $\Theta(u)=\left(\Theta^1(u),...,\Theta^m(u)\right)^{\top}$, Hiriart-Urruty and Lemar\'echal gives a dual optimal bound as follows. (See section 2.3 Chapter VII of~\cite{Slaters})
\begin{eqnarray}\label{eq:M0¡ª1}
\|p^*\|\leq\mu_0=\frac{(G+J)(\hat{u})-(G+J)(u^*)}{\min\limits_{1\leq j\leq m}\{-\Theta^j(\hat{u})\}},
\end{eqnarray}
where $\hat{u}$ is one Slater vector. In many cases, we can easy to get $\hat{u}$ and prior to estimate the lower bound $\underline{G+J}$ of $(G+J)(u^*)$, then we obtain one bound $\mu$ for the dual optimal set as following:
\begin{eqnarray}
\mu=\frac{(G+J)(\hat{u})-\underline{G+J}}{\min\limits_{1\leq j\leq m}\{-\Theta^j(\hat{u})\}}+1.
\end{eqnarray}
For constrained optimization problem Nedi\'c and Ozdaglar~\cite{bound1} recently use these bounds to design a dual algorithm that projects on a set containing the dual optimal solution. Shefi and Teboulle~\cite{bound2} also use these bounds to derive the rate of convergence results for the proximal method of multipliers for convex minimization.\\
\indent When $\mathbf{C}=\mathcal{K}_{\nu}^m$, Zhao and Zhu~\cite{ZhaoZhu2017} gives a dual optimal bound, and the following lemma shows that $\mu_0$ is computable.
\begin{lemma}\label{lemma:M}
If there exists a point $\hat{u}$ satisfying CQC condition for problem (P) and $\mathbf{C}=\mathcal{K}_{\nu}^{m+1}=\{x=(x_0,\overline{x})\in \RR\times\RR^{m}|x_0\geq\|\overline{x}\|_{\nu}\}$, then we have
\begin{eqnarray}\label{eq:M0}
\|p^*\|\leq \mu_0= m^{\max\{\frac{\omega-2}{2\omega},0\}}\cdot2^{\frac{1}{\omega}}\cdot\frac{(G+J)(\hat{u})-\underline{G+J}}{\theta_0-\|\overline{\theta}\|_\nu},
\end{eqnarray}
where $\frac{1}{\omega}+\frac{1}{\nu}=1$, $\underline{G+J}$ is the lower bound of $(G+J)(u^*)$ and $\Theta(\hat{u})=\left(\begin{array}{l}\theta_0\\\overline{\theta}\end{array}\right)$.
\end{lemma}
Then $$\mu=m^{\max\{\frac{\omega-2}{2\omega},0\}}\cdot2^{\frac{1}{\omega}}\cdot\frac{(G+J)(\hat{u})-\underline{G+J}}{\theta_0-\|\overline{\theta}\|_\nu}+1.$$
\indent Such bound plays a key role in our subsequence development. In particular, we use this bound to construct a new coordinate descent method by augmented Lagrangian to solve the nonlinear convex cone program. In order to treat the nonseparability of constraints, the dual update is need to compute a projection onto an intersection to the nonnegative orthant and a ball $\mathfrak{B}_\mu$ with radius $\mu$. However, this projection is very easy to implement as: $\mathcal{P}_\mu(y)=[\min(1,\mu/\|y\|)]y$.\\
\indent In the following sections, we propose a Stochastic Primal-Dual Coordinate method (SPDC) to solve NCCP, analyse the convergence and convergence rate of our method.
\section{Stochastic primal-dual coordinate method}\label{SPDC}
In this section, we propose a stochastic primal-dual coordinate descent algorithm to solve (P). Firstly, we introduce the core function $K(\cdot)$ satisfying the following assumption:
\begin{assumption}\label{assump2}
\rm{
$K$ is strongly convex with parameter $\beta$ and differentiable with its gradient Lipschitz continuous with parameter $B$ on $\mathbf{U}$.
}
\end{assumption}
\begin{remark}
Noted that $D(u,v)=K(u)-K(v)-\langle\nabla K(v),u-v\rangle$ is a Bregman like function (core function)~\cite{Mirror, CohenZ}. From Assumption~\ref{assump2} we have: $\frac{\beta}{2}\|u-v\|^2\leq D(u,v)\leq\frac{B}{2}\|u-v\|^2$.
\end{remark}
\indent Moreover, we assume that the sequence of parameter $\{\epsilon^k\}$ is non-increasing and satisfies:
\begin{equation}\label{para-choice-convergence}
0<\epsilon^{k+1}\leq\epsilon^{k}<\overline{\epsilon}<N\beta/\big{(}NB_{G}+\gamma\tau^2\big{)},\quad \sum_{k=0}^{+\infty}\epsilon^k=+\infty\quad\mbox{and}\quad\sum_{k=0}^{+\infty}(\epsilon^k)^2<+\infty.
\end{equation}
Then we introduce Stochastic Primal-Dual Coordinate Method (SPDC) for solving (P) in Algorithm~\ref{alg:SPDC}.
\begin{algorithm}
\caption{Stochastic Primal-Dual Coordinate Method (SPDC)}\label{alg:SPDC}
\begin{algorithmic}
\STATE{Initialize $u^0 \in \mathbf{U}$ and $p^0\in \mathbf{C^*}\cap\mathfrak{B}_{\mu}$}
\FOR{$k = 0,1,\cdots $,}
\STATE{Choose $i(k)$ from $\{1,2,\ldots,N\}$ with equal probability}
\STATE{Update}
\STATE{$u^{k+1}=\arg\min\limits_{u\in \mathbf{U}}\langle\nabla_{i(k)} G(u^{k}), u_{i(k)} \rangle + J_{i(k)}(u_{i(k)})+\langle
\Pi(p^k+\gamma\Theta(u^k)), \Theta_{i(k)}(u_{i(k)})\rangle+\frac{1}{\epsilon^k}D(u,u^k)$}
\STATE{Update $p^{k+1}=\mathcal{P}_\mu\bigg{(}\Pi\big{(}p^k+\gamma\Theta(u^{k+1})\big{)}\bigg{)}$}
\ENDFOR
\end{algorithmic}
\end{algorithm}
\newpage
In this algorithm, $\mathcal{P}_\mu$ denotes the projection on the ball $\mathfrak{B}_\mu$ which is computed easily. For the sake of brevity, let us set that $q^k=\Pi\big{(}p^k+\gamma\Theta(u^k)\big{)}$, $q^{k+1/2}=\Pi\big{(}p^k+\gamma\Theta(u^{k+1})\big{)}$.
Noted that the primal problem of algorithm can be expressed as
\begin{eqnarray}\label{eq:APk}
\qquad\mbox{(AP$^k$)}\;\min_{u\in \mathbf{U}}\langle\nabla_{i(k)} G(u^{k}), u_{i(k)}\rangle + J_{i(k)}(u_{i(k)})+ \langle q^k, \Theta_{i(k)}(u_{i(k)})\rangle+\frac{1}{\epsilon^k}D(u,u^k),
\end{eqnarray}
and the dual problem of algorithm is expressed as
\begin{equation}
p^{k+1}\leftarrow \mathcal{P}_\mu\bigg{(}q^{k+1/2}\bigg{)}.
\end{equation}
Note that if we choose an additive Bregman like function (or core function) respect to the space decomposition~\eqref{spacedecomposition} that is
$$K(u)=\sum_{i=1}^NK_i(u_i).$$
Then problem (AP$^k$) is just a small optimization problem for selected block $i(k)$. Specifically, taking $K(u)=\sum\limits_{i=1}\limits^{N}\frac{\|u_i\|^2}{2}$ for (AP$^k$), we perform only a block proximal gradient update for block $i(k)$, where we linearize the coupled function $G(u)$ and augmented Lagrangian term $\varphi(\Theta(u),p)$ and add the proximal term to it.\\
\indent In next sections, we will establish the convergence and convergence rate and probability complexity bounds of SPDC.
\section{Convergence analysis}\label{Convergence}
In this section, we will establish results about convergence of SPDC. Before proceeding, we first give the generalized equilibrium reformulation of saddle point formulation~\eqref{saddle point:L}:\\
Find $(u^*, p^*)\in\mathbf{U}\times\mathbf{C}^*$ such that
\begin{equation}\label{VIS}
\mbox{(EP):}\qquad L(u^*, p)-L(u, p^*)\leq 0, \forall u\in\mathbf{U}, p\in\mathbf{C}^*.
\end{equation}
Obviously, bifunction $L(u',p)-L(u,p')$ is convex in $u'$ and linear in $p'$ for given $u\in\mathbf{U}$, $p\in\mathbf{C}^*$.\\
\indent In algorithm SPDC, the indices $i(k)$, $k=0,1,2,\ldots$ are random variables. After $k$ iterations, SPDC method generates a random output $(u^{k+1}, p^{k+1})$. We denote by $\mathcal{F}_k$ is a filtration generated by the random variable $i(0),i(1),\ldots,i(k)$, i.e.,
$$\mathcal{F}_{k}\overset{def}{=}\{i(0),i(1),\ldots,i(k)\}, \mathcal{F}_{k}\subset\mathcal{F}_{k+1}.$$
\indent Additionaly, we define that $\mathcal{F}=(\mathcal{F}_{k})_{k\in\mathbb{N}}$,  $\mathbb{E}_{\mathcal{F}_{k+1}}=\mathbb{E}(\cdot|\mathcal{F}_{k})$ is the conditional expectation w.r.t. $\mathcal{F}_{k}$ and the conditional expectation in term of $i(k)$ given $i(0),i(1),\ldots,i(k-1)$ as $\mathbb{E}_{i(k)}$. \\
\indent Knowing $\mathcal{F}_{k-1}=\{i(0),i(1),\ldots,i(k-1)\}$, for given $u\in\mathbf{U}$, we have:
\begin{eqnarray}
&&\qquad\qquad\mathbb{E}_{i(k)}\langle\nabla_{i(k)}G(u^{k}),(u^{k}-u)_{i(k)}\rangle=\frac{1}{N}\langle\nabla G(u^{k}),u^{k}-u\rangle\geq\frac{1}{N}\big{[}G(u^k)-G(u)\big{]};\label{expectationG}\\ &&\qquad\qquad\mathbb{E}_{i(k)}\big{[}J_{i(k)}(u_{i(k)}^{k})-J_{i(k)}(u_{i(k)})\big{]}=\frac{1}{N}\big{[}J(u^{k})-J(u)\big{]};\label{expectationJ}\\ &&\qquad\qquad\mathbb{E}_{i(k)}\langle q^k,\Theta_{i(k)}(u_{i(k)}^{k})-\Theta_{i(k)}(u_{i(k)})\rangle=\frac{1}{N}\langle q^k,\Theta(u^{k})-\Theta(u)\rangle.\label{expectationP}
\end{eqnarray}
For the sequence $\{(u^k,p^k)\}$ generated by algorithm SPDC, the following function provides one upper bound for $\|u^k-u^*\|^2$:
\begin{eqnarray}\label{Lambda:bound}
\Lambda^k(u^{k},p^k)&=&D(u^*,u^k)+\frac{\epsilon^k}{2\gamma N}\|p^*-p^k\|^2+\epsilon^k\big{(}L(u^k,p^*)-L(u^*,p^*)\big{)}\\
                    &\geq&\frac{\beta}{2}\|u^k-u^*\|^2+\frac{\epsilon^k}{2\gamma N}\|p^*-p^k\|^2\nonumber\\
                    &\geq&\frac{\beta}{2}\|u^k-u^*\|^2.\nonumber
\end{eqnarray}
The boundness of the sequence $\{(u^k,p^k)\}$ generated by SPDC will play an important role. To obtain this boundness, we need the following lemma.
\begin{lemma}\label{lemma:bound1} {\bf (Descent inequalities for bifunction values with Lagrangian)}
Let Assumption~\ref{assump1} and~\ref{assump2} hold, $\{(u^k,p^k)\}$ is generated by SPDC, the parameter sequence $\{\epsilon^k\}$  satisfies condition~\eqref{para-choice-convergence}. Then $\forall p\in\mathbf{C}^*\cap\mathfrak{B}_{\mu}$ it holds that
\begin{eqnarray*}
\mbox{{\rm(i)}}\quad\frac{\epsilon^k}{N}\big{[}L(u^{k}, q^{k})-L(u^*, q^{k})\big{]}&\leq&\big{[}D(u^*,u^k)+\epsilon^k\big{(}L(u^{k},p^*)-L(u^*,p^*)\big{)}\big{]}\\
&&-\mathbb{E}_{i(k)}\big{[}D(u^*,u^{k+1})+\epsilon^k\big{(}L(u^{k+1},p^*)-L(u^*,p^*)\big{)}\big{]}\\
&&+\epsilon^k\mathbb{E}_{i(k)}\langle q^k-p^*,\Theta(u^{k})-\Theta(u^{k+1})\rangle-\frac{\beta-\epsilon^kB_G}{2}\mathbb{E}_{i(k)}\|u^k-u^{k+1}\|^2\\
\mbox{{\rm(ii)}}\quad\frac{\epsilon^k}{N}\big{[}L(u^k,p)-L(u^{k},q^{k})\big{]}&\leq&\frac{1}{2\gamma N}\big{[}\epsilon^k\|p-p^k\|^2-\epsilon^{k+1}\|p-p^{k+1}\|^2\big{]}\\
&&-\frac{\epsilon^k}{N}\langle q^k-p,\Theta(u^{k})-\Theta(u^{k+1})\rangle+\frac{\gamma\tau^2\epsilon^k}{2N}\|u^{k}-u^{k+1}\|^2-\frac{\epsilon^k}{2\gamma N}\|q^k-p^{k}\|^2\\
\mbox{{\rm(iii)}}\;\;\frac{\epsilon^k}{N}\big{[}L(u^{k}, p^{*})-L(u^*, q^{k})\big{]}&\leq&(1+\eta_1(\epsilon^k)^2)\Lambda^k(u^{k},p^k)-\mathbb{E}_{i(k)}\Lambda^{k+1}(u^{k+1},p^{k+1})\\
&&+\eta_2(\epsilon^k)^2-\eta_3\mathbb{E}_{i(k)}\|u^k-u^{k+1}\|^2-\frac{\epsilon^k}{2\gamma N}\|q^k-p^{k}\|^2
\end{eqnarray*}
where $\eta_1=\frac{\gamma\tau^3(N-1)^2(1+\gamma\tau)}{N\beta[N\beta-\bar{\epsilon}(NB_G+\gamma\tau^2)]}$, $\eta_2=\frac{\tau^2(N-1)^2(1+\gamma\tau)\mu^2}{N[N\beta-\bar{\epsilon}(NB_G+\gamma\tau^2)]}$ and $\eta_3=\frac{(\bar{\epsilon}-\epsilon^0)(NB_G+\gamma\tau^2)}{2N}$.
\end{lemma}
\begin{proof}
(i) Firstly, for all $u\in\mathbf{U}$, the unique solution $u^{k+1}$ of the primal problem~\eqref{eq:APk} is characterized by the following variational inequality:
\begin{eqnarray}\label{eq:primal_VI}
\\
\langle\nabla_{i(k)}G(u^{k}),(u^{k+1}-u)_{i(k)}\rangle+J_{i(k)}(u_{i(k)}^{k+1})-J_{i(k)}(u_{i(k)})+\langle q^k,\Theta_{i(k)}(u_{i(k)}^{k+1})-\Theta_{i(k)}(u_{i(k)})\rangle\nonumber\\
+\frac{1}{\epsilon^k}\langle\nabla K(u^{k+1})-\nabla K(u^k), u^{k+1}-u\rangle\leq 0,\nonumber
\end{eqnarray}
which follows that
\begin{eqnarray}\label{eq:primal_VI2}
\langle\nabla_{i(k)}G(u^{k}),\big{(}u^{k}-u-(u^{k}-u^{k+1})\big{)}_{i(k)}\rangle+J_{i(k)}(u_{i(k)}^{k})-J_{i(k)}(u_{i(k)})\\
-\big{(}J_{i(k)}(u_{i(k)}^{k})-J_{i(k)}(u_{i(k)}^{k+1})\big{)}\nonumber\\
+\langle q^k,\Theta_{i(k)}(u_{i(k)}^{k})-\Theta_{i(k)}(u_{i(k)})-\big{(}\Theta_{i(k)}(u_{i(k)}^{k})-\Theta_{i(k)}(u_{i(k)}^{k+1})\big{)}\rangle\nonumber\\
+\frac{1}{\epsilon^k}\langle\nabla K(u^{k+1})-\nabla K(u^k), u^{k+1}-u\rangle\leq 0.\nonumber
\end{eqnarray}
Observing that $\langle\nabla_{i(k)}G(u^{k}),(u^{k}-u^{k+1})_{i(k)}\rangle=\langle\nabla G(u^{k}),u^{k}-u^{k+1}\rangle$, $J_{i(k)}(u_{i(k)}^{k})-J_{i(k)}(u_{i(k)}^{k+1})=J(u^{k})-J(u^{k+1})$ and $\Theta_{i(k)}(u_{i(k)}^{k})-\Theta_{i(k)}(u_{i(k)}^{k+1})=\Theta(u^{k})-\Theta(u^{k+1})$, from~\eqref{eq:primal_VI2}, we have that
\begin{eqnarray}\label{eq:primal_bound1}
\\
\langle\nabla_{i(k)}G(u^{k}),(u^{k}-u)_{i(k)}\rangle+J_{i(k)}(u_{i(k)}^{k})-J_{i(k)}(u_{i(k)})+\langle q^k,\Theta_{i(k)}(u_{i(k)}^{k})-\Theta_{i(k)}(u_{i(k)})\rangle\nonumber\\
\leq\langle\nabla G(u^{k}),u^{k}-u^{k+1}\rangle+J(u^{k})-J(u^{k+1})+\langle q^k,\Theta(u^{k})-\Theta(u^{k+1})\rangle\nonumber\\
+\frac{1}{\epsilon^k}\langle\nabla K(u^{k+1})-\nabla K(u^k), u-u^{k+1}\rangle.\nonumber
\end{eqnarray}
By Lemma~\ref{lemma:fphi}, we have that
\begin{eqnarray}\label{eq:primal_bound1-1}
\\
\langle\nabla_{i(k)}G(u^{k}),(u^{k}-u)_{i(k)}\rangle+J_{i(k)}(u_{i(k)}^{k})-J_{i(k)}(u_{i(k)})+\langle q^k,\Theta_{i(k)}(u_{i(k)}^{k})-\Theta_{i(k)}(u_{i(k)})\rangle\nonumber\\
\leq(G+J)(u^k)-(G+J)(u^{k+1})+\frac{B_G}{2}\|u^{k}-u^{k+1}\|^2+\langle q^k,\Theta(u^{k})-\Theta(u^{k+1})\rangle\nonumber\\
+\frac{1}{\epsilon^k}\langle\nabla K(u^{k+1})-\nabla K(u^k), u-u^{k+1}\rangle.\nonumber
\end{eqnarray}
The simple algebraic operation follows that
\begin{eqnarray}\label{eq:primal_bound2}
\\
\frac{1}{\epsilon^k}\langle \nabla K(u^{k+1})-\nabla K(u^k),u-u^{k+1}\rangle
&=&\frac{1}{\epsilon^k}\big{[}D(u,u^k)-D(u,u^{k+1})-D(u^{k+1},u^k)\big{]}\nonumber\\
&\leq&\frac{1}{\epsilon^k}\big{[}D(u,u^k)-D(u,u^{k+1})\big{]}-\frac{\beta}{2\epsilon^k}\|u^k-u^{k+1}\|^2.\nonumber
\end{eqnarray}
Combining~\eqref{eq:primal_bound1-1} and~\eqref{eq:primal_bound2} with $u=u^*$, we obtain that
\begin{eqnarray}\label{eq:primal_bound3}
\\
\langle\nabla_{i(k)}G(u^{k}),(u^{k}-u^*)_{i(k)}\rangle+J_{i(k)}(u_{i(k)}^{k})-J_{i(k)}(u_{i(k)}^*)+\langle q^k,\Theta_{i(k)}(u_{i(k)}^{k})-\Theta_{i(k)}(u_{i(k)}^*)\rangle\nonumber\\
\leq\frac{1}{\epsilon^k}\big{[}D(u^*,u^k)-D(u^*,u^{k+1})\big{]}+(G+J)(u^k)-(G+J)(u^{k+1})\nonumber\\
+\langle q^k,\Theta(u^{k})-\Theta(u^{k+1})\rangle-\frac{\beta-\epsilon^kB_G}{2\epsilon^k}\|u^k-u^{k+1}\|^2.\nonumber
\end{eqnarray}
Take expectation with respect to $i(k)$ on both side of~\eqref{eq:primal_bound3}, together the conditional expectation~\eqref{expectationG}-\eqref{expectationP}, we get
\begin{eqnarray}\label{eq:primal_bound4}
&&\frac{1}{N}\big{[}(G+J)(u^{k})-(G+J)(u^*)+\langle q^k,\Theta(u^{k})-\Theta(u^*)\rangle\big{]}\\
&\leq&\frac{1}{\epsilon^k}\big{[}D(u^*,u^k)-\mathbb{E}_{i(k)}D(u^*,u^{k+1})\big{]}+(G+J)(u^k)-\mathbb{E}_{i(k)}(G+J)(u^{k+1})\nonumber\\
&&+\mathbb{E}_{i(k)}\langle q^k,\Theta(u^{k})-\Theta(u^{k+1})\rangle-\frac{\beta-\epsilon^kB_G}{2\epsilon^k}\mathbb{E}_{i(k)}\|u^k-u^{k+1}\|^2.\nonumber
\end{eqnarray}
Furthermore, we have that
\begin{eqnarray*}\label{eq:primal_bound5}
&&\frac{\epsilon^k}{N}\big{[}L(u^{k}, q^{k})-L(u^*, q^{k})\big{]}\nonumber\\
&=&\frac{\epsilon^k}{N}\big{[}(G+J)(u^{k})-(G+J)(u^*)+\langle q^k,\Theta(u^{k})-\Theta(u^*)\rangle\big{]}\nonumber\\
&\leq&\big{[}D(u^*,u^k)+\epsilon^k\big{(}L(u^{k},p^*)-L(u^*,p^*)\big{)}\big{]}-\mathbb{E}_{i(k)}\big{[}D(u^*,u^{k+1})+\epsilon^k\big{(}L(u^{k+1},p^*)-L(u^*,p^*)\big{)}\big{]}\nonumber\qquad\\
&&+\epsilon^k\mathbb{E}_{i(k)}\langle q^k-p^*,\Theta(u^{k})-\Theta(u^{k+1})\rangle-\frac{\beta-\epsilon^kB_G}{2}\mathbb{E}_{i(k)}\|u^k-u^{k+1}\|^2.\qquad\mbox{(from~\eqref{eq:primal_bound4})}\nonumber
\end{eqnarray*}
Then statement (i) is provided.\\
(ii) By simple operations the following equality holds for all $p\in\mathbf{C}^*$:
\begin{eqnarray}\label{eq:dual_bound1}
\qquad\langle p-q^k,\Theta(u^{k+1})\rangle
&=&\frac{1}{\gamma}\langle p-q^{k+1/2}, p^k+\gamma \Theta(u^{k+1})-q^{k+1/2}\rangle\\
&&+\frac{1}{\gamma}\langle p-q^{k+1/2}, q^{k+1/2}-p^k\rangle+\langle q^{k+1/2}-q^k,\Theta(u^{k+1})\rangle.\nonumber
\end{eqnarray}
From the property of projection~\eqref{eq:Projecproperty3} with $x=p^k+\gamma \Theta(u^{k+1})$ and $y=p$, the first term of right hand side of~\eqref{eq:dual_bound1} is negative:
\begin{equation}\label{eq:dual_bound2}
\frac{1}{\gamma}\langle p-q^{k+1/2}, p^k+\gamma \Theta(u^{k+1})-q^{k+1/2}\rangle\leq 0.
\end{equation}
Second, using Proposition~\ref{proposition} with $x=\gamma\Theta(u^{k+1})$, $y=\gamma\Theta(u^k)$, $z=p^k$, we obtain one upper bound of the third term of~\eqref{eq:dual_bound1}:
\begin{eqnarray}\label{eq:dual_bound3}
\\
\langle q^{k+1/2}-q^k,\Theta(u^{k+1})\rangle
&\leq&\frac{1}{2\gamma}[\|q^{k+1/2}-p^k\|^2-\|q^k-p^k\|^2]+\frac{\gamma}{2}\|\Theta(u^k)-\Theta(u^{k+1})\|^2\nonumber\\
&\leq&\frac{1}{2\gamma}[\|q^{k+1/2}-p^k\|^2-\|q^k-p^k\|^2]+\frac{\gamma\tau^2}{2}\|u^k-u^{k+1}\|^2.\nonumber
\end{eqnarray}
Thus from~\eqref{eq:dual_bound1}, for any $p\in\mathbf{C}^*\cap \mathfrak{B}_\mu$, we have
\begin{eqnarray}\label{eq:dual_bound4}
\\
&&\langle p-q^k,\Theta(u^{k+1})\rangle\nonumber\\
&\leq&\frac{1}{\gamma}\langle p-q^{k+1/2}, q^{k+1/2}-p^k\rangle+\frac{1}{2\gamma}[\|q^{k+1/2}-p^k\|^2-\|q^k-p^k\|^2]+\frac{\gamma\tau^2}{2}\|u^k-u^{k+1}\|^2\nonumber\\
&\leq&\frac{1}{2\gamma}[\|p-p^k\|^2-\|p-q^{k+1/2}\|^2]-\frac{1}{2\gamma}\|q^k-p^k\|^2+\frac{\gamma\tau^2}{2}\|u^k-u^{k+1}\|^2\nonumber\\
                            &\leq&\frac{1}{2\gamma}[\|p-p^k\|^2-\|p-p^{k+1}\|^2]-\frac{1}{2\gamma}\|q^k-p^k\|^2+\frac{\gamma\tau^2}{2}\|u^k-u^{k+1}\|^2.\nonumber\\
                              &&\qquad\qquad\qquad\mbox{(since $p^{k+1}=\mathcal{P}_\mu(q^{k+1/2})$ and by property~\eqref{eq:Projecproperty4} of projection)}\nonumber
\end{eqnarray}
Therefore, we have
\begin{eqnarray}\label{eq:dual_bound4-1}
&&\langle p-q^k,\Theta(u^{k})\rangle\\
&=&\langle p-q^k,\Theta(u^{k+1})\rangle+\langle p-q^k,\Theta(u^{k})-\Theta(u^{k+1})\rangle\nonumber\\
&\leq&\frac{1}{2\gamma}\|p-p^k\|^2-\frac{1}{2\gamma}\|p-p^{k+1}\|^2-\langle q^k-p, \Theta(u^{k})-\Theta(u^{k+1})\rangle\nonumber\\
&&+\frac{\gamma\tau^2}{2}\|u^{k}-u^{k+1}\|^2-\frac{1}{2\gamma}\|q^k-p^k\|^2\nonumber
\end{eqnarray}
Then,~\eqref{eq:dual_bound4-1} yields that
\begin{eqnarray}\label{eq:dual_bound5}
&&\frac{\epsilon^k}{N}\big{[}L(u^k,p)-L(u^{k},q^{k})\big{]}\\
&=&\frac{\epsilon^k}{N}\langle p-q^k,\Theta(u^{k})\rangle\nonumber\\
&\leq&\frac{\epsilon^k}{2\gamma N}\|p-p^k\|^2-\frac{\epsilon^k}{2\gamma N}\|p-p^{k+1}\|^2-\frac{\epsilon^k}{N}\langle q^k-p, \Theta(u^{k})-\Theta(u^{k+1})\rangle\nonumber\\
&&+\frac{\gamma\tau^2\epsilon^k}{2N}\|u^{k}-u^{k+1}\|^2-\frac{\epsilon^k}{2\gamma N}\|q^k-p^k\|^2\nonumber\\
&\leq&\frac{1}{2\gamma N}\big{(}\epsilon^k\|p-p^k\|^2-\epsilon^{k+1}\|p-p^{k+1}\|^2\big{)}-\frac{\epsilon^k}{N}\langle q^k-p, \Theta(u^{k})-\Theta(u^{k+1})\rangle\nonumber\\
&&+\frac{\gamma\tau^2\epsilon^k}{2N}\|u^{k}-u^{k+1}\|^2-\frac{\epsilon^k}{2\gamma N}\|q^k-p^k\|^2\qquad\mbox{(since $\epsilon^{k+1}\leq\epsilon^k$)}.\nonumber
\end{eqnarray}
Then statement (ii) is provided.\\
(iii) First, we observe that
\begin{eqnarray}\label{eq:all_bound1}
\qquad&&\langle q^k-p^*,\Theta(u^{k})-\Theta(u^{k+1})\rangle\\
&=&\langle\Pi\big{(}p^k+\gamma\Theta(u^k)\big{)}-\Pi\big{(}p^*+\gamma\Theta(u^*)\big{)},\Theta(u^{k})-\Theta(u^{k+1})\rangle\nonumber\\
&\leq&\|\Pi\big{(}p^k+\gamma\Theta(u^k)\big{)}-\Pi\big{(}p^*+\gamma\Theta(u^*)\big{)}\|\cdot\|\Theta(u^{k})-\Theta(u^{k+1})\|\nonumber\\
&\leq&\tau\big{[}\|p^*-p^k\|\cdot\|u^{k}-u^{k+1}\|+\gamma\tau\|u^*-u^k\|\cdot\|u^{k}-u^{k+1}\|\big{]}\nonumber\\
&\leq&\frac{\tau}{2}\left[\kappa\epsilon^k\big{(}\|p^*-p^k\|^2+\gamma\tau\|u^*-u^k\|^2\big{)}+\frac{1+\gamma\tau}{\kappa\epsilon^k}\|u^{k}-u^{k+1}\|^2\right]\nonumber\\
&&\qquad\qquad\qquad\qquad\quad\mbox{(from the H\"older's inequality $xy\leq(\kappa x^2+\frac{y^2}{\kappa})/2$, $\kappa>0$)}\nonumber\\
&\leq&\frac{\tau}{2}\left[\kappa\epsilon^k\big{(}2\mu^2+\gamma\tau\|u^*-u^k\|^2\big{)}+\frac{1+\gamma\tau}{\kappa\epsilon^k}\|u^{k}-u^{k+1}\|^2\right]\;\;\mbox{(since $p^*,p^k\in\mathbf{C}^*\cap\mathfrak{B}_{\mu}$)}\nonumber
\end{eqnarray}
For statement (ii), let $p=p^*$, and take the expectation of $i(k)$ on both sides, we have
\begin{eqnarray}\label{eq:all_bound2}
\\
&&\frac{\epsilon^k}{N}\big{[}L(u^k,p^*)-L(u^{k},q^{k})\big{]}\nonumber\\
&\leq&\frac{1}{2\gamma N}\big{(}\epsilon^k\|p^*-p^k\|^2-\epsilon^{k+1}\mathbb{E}_{i(k)}\|p^*-p^{k+1}\|^2\big{)}-\frac{\epsilon^k}{N}\mathbb{E}_{i(k)}\langle q^k-p^*, \Theta(u^{k})-\Theta(u^{k+1})\rangle\nonumber\\
&&+\frac{\gamma\tau^2\epsilon^k}{2N}\mathbb{E}_{i(k)}\|u^{k}-u^{k+1}\|^2-\frac{\epsilon^k}{2\gamma N}\|q^k-p^k\|^2.\nonumber
\end{eqnarray}
Hence, summing statement (i) and~\eqref{eq:all_bound2}, it follows
\begin{eqnarray*}
&&\frac{\epsilon^k}{N}\big{[}L(u^k,p^*)-L(u^{*},q^{k})\big{]}\\
&\leq&\big{[}D(u^*,u^k)+\epsilon^k\big{(}L(u^{k},p^*)-L(u^*,p^*)\big{)}+\frac{\epsilon^k}{2\gamma N}\|p^*-p^k\|^2\big{]}\\
&&-\mathbb{E}_{i(k)}\big{[}D(u^*,u^{k+1})+\epsilon^k\big{(}L(u^{k+1},p^*)-L(u^*,p^*)\big{)}+\frac{\epsilon^{k+1}}{2\gamma N}\|p^*-p^{k+1}\|^2\big{]}\nonumber\\
&&+\mathbb{E}_{i(k)}\big{[}\frac{(N-1)\epsilon^k}{N}\langle q^k-p^*,\Theta(u^{k})-\Theta(u^{k+1})\rangle-\frac{N\beta-\epsilon^k(NB_G+\gamma\tau^2)}{2N}\|u^k-u^{k+1}\|^2\big{]}\nonumber\\
&&-\frac{\epsilon^k}{2\gamma N}\|q^k-p^{k}\|^2\nonumber\\
&\leq&\Lambda^k(u^{k},p^k)-\mathbb{E}_{i(k)}\Lambda^{k+1}(u^{k+1},p^{k+1})\\
&&+\mathbb{E}_{i(k)}\big{[}\frac{(N-1)\epsilon^k}{N}\langle q^k-p^*,\Theta(u^{k})-\Theta(u^{k+1})\rangle-\frac{N\beta-\epsilon^k(NB_G+\gamma\tau^2)}{2N}\|u^k-u^{k+1}\|^2\big{]}\nonumber\\
&&-\frac{\epsilon^k}{2\gamma N}\|q^k-p^{k}\|^2\qquad\qquad\qquad\qquad\qquad\qquad\qquad\qquad\qquad\qquad\qquad\mbox{(since $\epsilon^{k+1}\leq\epsilon^k$)}.
\end{eqnarray*}
By~\eqref{eq:all_bound1}, we have that
\begin{eqnarray*}
&&\frac{\epsilon^k}{N}\big{[}L(u^k,p^*)-L(u^{*},q^{k})\big{]}\nonumber\\
&\leq&\Lambda^k(u^{k},p^k)-\mathbb{E}_{i(k)}\Lambda^{k+1}(u^{k+1},p^{k+1})+\mathbb{E}_{i(k)}\big{[}\frac{(N-1)\epsilon^k}{N}\langle q^k-p^*,\Theta(u^{k})-\Theta(u^{k+1})\rangle\nonumber\\
&&-\frac{N\beta-\epsilon^k(NB_G+\gamma\tau^2)}{2N}\|u^k-u^{k+1}\|^2\big{]}-\frac{\epsilon^k}{2\gamma N}\|q^k-p^{k}\|^2\nonumber\\
&\leq&\Lambda^k(u^{k},p^k)-\mathbb{E}_{i(k)}\Lambda^{k+1}(u^{k+1},p^{k+1})+\frac{\kappa\tau(N-1)(\epsilon^k)^2}{2N}\big{[}2\mu^2+\gamma\tau\|u^*-u^k\|^2\big{]}\nonumber\\
&&+\frac{\tau(1+\gamma\tau)(N-1)-\kappa\big{(} N\beta-\epsilon^k(NB_G+\gamma\tau^2)\big{)}}{2\kappa N}\mathbb{E}_{i(k)}\|u^k-u^{k+1}\|^2-\frac{\epsilon^k}{2\gamma N}\|q^k-p^{k}\|^2.
\end{eqnarray*}
Taking $\kappa=\frac{\tau(1+\gamma\tau)(N-1)}{N\beta-\bar{\epsilon}(NB_G+\gamma\tau^2)}$, by definition of $\eta_1$, $\eta_2$ and $\eta_3$ of this statement, we have that
\begin{eqnarray*}
\frac{\epsilon^k}{N}\big{[}L(u^k,p^*)-L(u^{*},q^{k})\big{]}
&\leq&\Lambda^k(u^{k},p^k)-\mathbb{E}_{i(k)}\Lambda^{k+1}(u^{k+1},p^{k+1})+\eta_1(\epsilon^k)^2\frac{\beta}{2}\|u^*-u^k\|^2+\eta_2(\epsilon^k)^2\nonumber\\
&&-\eta_3\mathbb{E}_{i(k)}\|u^k-u^{k+1}\|^2-\frac{\epsilon^k}{2\gamma N}\|q^k-p^{k}\|^2.
\end{eqnarray*}
By~\eqref{Lambda:bound}, we have $\frac{\beta}{2}\|u^*-u^k\|^2\leq\Lambda^k(u^{k},p^k)$. Then we obtain statement (iii).
\end{proof}
\indent Based Lemma~\ref{lemma:bound1}, we establish the following proposition for the boundness of sequence $\{u^k\}$ generated by SPDC.\\
\begin{proposition}[{\bf Almost surely boundness for $\{u^k\}$}]\label{proposition2}
Let assumptions of Lemma~\ref{lemma:bound1} hold, then
\begin{itemize}
\item[{\rm(i)}] $\sum\limits_{k=0}\limits^{+\infty}\mathbb{E}_{i(k)}\|u^k-u^{k+1}\|^2<+\infty$ $\mbox{a.s.}$ and $\sum\limits_{k=0}\limits^{+\infty}\epsilon^k\|q^k-p^{k}\|^2<+\infty$ $\mbox{a.s.}$;
\item[{\rm(ii)}] The sequence $\{u^{k}\}$ generated by SPDC is almost surely bounded.
\end{itemize}
\end{proposition}
\begin{proof}
(i) From statement (iii) of Lemma~\ref{lemma:bound1}, we have
\begin{eqnarray}\label{eq:bound2_1_8}
\mathbb{E}_{i(k)}\Lambda^{k+1}(u^{k+1},p^{k+1})\leq(1+\eta_1(\epsilon^k)^2)\Lambda^k(u^{k},p^k)
+\eta_2(\epsilon^k)^2-S_k,
\end{eqnarray}
where $S_k=\frac{\epsilon^k}{N}\big{[}L(u^k,p^*)-L(u^{*},q^{k})\big{]}+\eta_3\mathbb{E}_{i(k)}\|u^k-u^{k+1}\|^2+\frac{\epsilon^k}{2\gamma N}\|q^k-p^{k}\|^2$ is positive.\\
By the Robbins-Siegmund Lemma~\cite{RS}, we obtain that $\lim\limits_{k\rightarrow+\infty}\Lambda^k(u^k,p^k)$ almost surely exists, $\sum\limits_{k=0}\limits^{+\infty}\mathbb{E}_{i(k)}\|u^k-u^{k+1}\|^2<+\infty$ $\mbox{a.s.}$ and $\sum\limits_{k=0}\limits^{+\infty}\epsilon^k\|q^k-p^{k}\|^2<+\infty$ $\mbox{a.s.}$.\\
(ii) Since $\lim\limits_{k\rightarrow+\infty}\Lambda^k(u^{k},p^k)$ almost surely exists, then $\Lambda^k(u^{k},p^k)$ is almost surely bounded. From~\eqref{Lambda:bound}, we have $\Lambda^k(u^{k},p^k)\geq\frac{\beta}{2}\|u^k-u^*\|^2$,
which implies the sequence $\{u^k\}$ is almost surely bounded.
\end{proof}
\indent Before the study of the convergence of SPDC, we start with following lemma about augmented Lagrangian $L_\gamma$, which will be used in the forthcoming convergence analysis.
\begin{lemma}\label{lemma:bound2} {\bf(Descent inequalities for the bifunction values with augmented Lagrangian)}
Suppose assumptions of Lemma~\ref{lemma:bound1} hold, then we have the following assertions:
\begin{eqnarray*}
\mbox{{\rm(i)}}\quad\frac{\epsilon^k}{N}\big{[}L_\gamma(u^{k}, p^{k})-L_\gamma(u^*, p^{k})\big{]}&\leq&\big{[}\Lambda^k(u^{k},p^k)-\mathbb{E}_{i(k)}\Lambda^{k+1}(u^{k+1},p^{k+1})\big{]}-\frac{1}{2\gamma N}\big{[}\epsilon^k\|p^*-p^k\|^2-\epsilon^{k+1}\mathbb{E}_{i(k)}\|p^*-p^{k+1}\|^2\big{]}\\
&&+\epsilon^k\tau\|q^k-p^*\|\cdot\mathbb{E}_{i(k)}\|u^{k}-u^{k+1}\|\\
\mbox{{\rm(ii)}}\quad\frac{\epsilon^k}{N}\big{[}L_\gamma(u^k,p^*)-L_\gamma(u^{k},p^{k})\big{]}&\leq&\frac{1}{2\gamma N}\big{[}\epsilon^k\|p^*-p^k\|^2-\epsilon^{k+1}\mathbb{E}_{i(k)}\|p^*-p^{k+1}\|^2\big{]}+\frac{\epsilon^k}{\gamma N}\|q^k-p^k\|^2\\
&&+\frac{\epsilon^0\gamma\tau^2}{N}\mathbb{E}_{i(k)}\|u^{k}-u^{k+1}\|^2+\frac{\epsilon^k\tau}{N}\|p^*-p^k\|\cdot\mathbb{E}_{i(k)}\|u^k-u^{k+1}\|\\
\mbox{{\rm(iii)}}\quad\frac{\epsilon^k}{N}\big{[}L_\gamma(u^{k}, p^{*})-L_\gamma(u^*, p^{*})\big{]}&\leq&\Lambda^k(u^{k},p^k)-\mathbb{E}_{i(k)}\Lambda^{k+1}(u^{k+1},p^{k+1})+\frac{\epsilon^k}{\gamma N}\|q^k-p^k\|^2+\frac{\epsilon^0\gamma\tau^2}{N}\mathbb{E}_{i(k)}\|u^{k}-u^{k+1}\|^2\\&&+h_1(u^k,p^k)(\epsilon^k)^2
\end{eqnarray*}
where $h_1(u^k,p^k)=\frac{2\tau}{\beta}\left(\|q^k-p^*\|+\frac{1}{N}\|p^*-p^k\|\right)\left(\|\nabla G(u^{k})\|+\|r^k\|+\tau\| q^k\|\right)$\\ and $r^k\in\partial J(u^k)$.
\end{lemma}
\begin{proof}
(i) By the convexity of $\varphi(\theta,p)$ in $\theta$, we have
\begin{eqnarray}\label{eq:AL_primal_1}
&&\frac{\epsilon^k}{N}\big{[}L_\gamma(u^{k},p^{k})-L_\gamma(u^*,p^{k})\big{]}\\
&=&\frac{\epsilon^k}{N}\big{[}(G+J)(u^k)-(G+J)(u^*)+\varphi(\Theta(u^k),p^k)-\varphi(\Theta(u^*),p^k)\big{]}\nonumber\\
                                                                           &\leq&\frac{\epsilon^k}{N}\big{[}(G+J)(u^k)-(G+J)(u^*)+\langle q^k,\Theta(u^k)-\Theta(u^*)\rangle\big{]}\nonumber\\
                                                                           &&\qquad\qquad\qquad\qquad\qquad\qquad\qquad\qquad\mbox{(since $\nabla\varphi_{\theta}\big{(}\Theta(u^k),p^k\big{)}=q^k$)}\nonumber\\
                                                                           &=&\frac{\epsilon^k}{N}\big{[}L(u^{k}, q^{k})-L(u^*, q^{k})\big{]}.\nonumber
\end{eqnarray}
Hence, from (i) of Lemma~\ref{lemma:bound1}, we get
\begin{eqnarray*}\label{eq:AL_primal_2}
&&\frac{\epsilon^k}{N}\big{[}L_\gamma(u^{k},p^{k})-L_\gamma(u^*,p^{k})\big{]}\\
&\leq&\big{[}D(u^*,u^k)+\epsilon^k\big{(}L(u^{k},p^*)-L(u^*,p^*)\big{)}\big{]}-\mathbb{E}_{i(k)}\big{[}D(u^*,u^{k+1})+\epsilon^k\big{(}L(u^{k+1},p^*)-L(u^*,p^*)\big{)}\big{]}\\
&&+\epsilon^k\mathbb{E}_{i(k)}\langle q^k-p^*,\Theta(u^{k})-\Theta(u^{k+1})\rangle-\frac{\beta-\epsilon^kB_G}{2}\mathbb{E}_{i(k)}\|u^k-u^{k+1}\|^2\\
&\leq&\big{[}\Lambda^k(u^{k},p^k)-\mathbb{E}_{i(k)}\Lambda^{k+1}(u^{k+1},p^{k+1})\big{]}-\frac{1}{2\gamma N}\big{[}\epsilon^k\|p^*-p^k\|^2-\epsilon^{k+1}\mathbb{E}_{i(k)}\|p^*-p^{k+1}\|^2\big{]}\nonumber\\
&&+\epsilon^k\mathbb{E}_{i(k)}\langle q^k-p^*,\Theta(u^{k})-\Theta(u^{k+1})\rangle-\frac{\beta-\epsilon^kB_G}{2}\mathbb{E}_{i(k)}\|u^k-u^{k+1}\|^2\qquad\mbox{(since $\epsilon^{k+1}\leq\epsilon^k$)}\nonumber\\
&\leq&\big{[}\Lambda^k(u^{k},p^k)-\mathbb{E}_{i(k)}\Lambda^{k+1}(u^{k+1},p^{k+1})\big{]}-\frac{1}{2\gamma N}\big{[}\epsilon^k\|p^*-p^k\|^2-\epsilon^{k+1}\mathbb{E}_{i(k)}\|p^*-p^{k+1}\|^2\big{]}\nonumber\\
&&+\epsilon^k\tau\|q^k-p^*\|\cdot\mathbb{E}_{i(k)}\|u^{k}-u^{k+1}\|-\frac{\beta-\epsilon^kB_G}{2}\mathbb{E}_{i(k)}\|u^k-u^{k+1}\|^2\nonumber\\
&\leq&\big{[}\Lambda^k(u^{k},p^k)-\mathbb{E}_{i(k)}\Lambda^{k+1}(u^{k+1},p^{k+1})\big{]}-\frac{1}{2\gamma N}\big{[}\epsilon^k\|p^*-p^k\|^2-\epsilon^{k+1}\mathbb{E}_{i(k)}\|p^*-p^{k+1}\|^2\big{]}\nonumber\\
&&+\epsilon^k\tau\|q^k-p^*\|\cdot\mathbb{E}_{i(k)}\|u^{k}-u^{k+1}\|.\nonumber
\end{eqnarray*}
(ii) The concavity of $\varphi(\theta,p)$ in $p$ yields that
\begin{eqnarray}\label{eq:AL_dual_1}
\\
\frac{\epsilon^k}{N}\big{[}L_\gamma(u^{k},p^{*})-L_\gamma(u^k,p^{k})\big{]}&=&\frac{\epsilon^k}{N}\bigg{(}\varphi\big{(}\Theta(u^k),p^*\big{)}-\varphi\big{(}\Theta(u^k),p^k\big{)}\bigg{)}\nonumber\\
&\leq&\frac{\epsilon^k}{\gamma N}\langle q^k-p^k,p^*-p^k\rangle\;\;\mbox{(since $\nabla\varphi_p\big{(}\Theta(u^k),p^k\big{)}=\frac{q^k-p^k}{\gamma}$)}\nonumber\\
&=&\frac{\epsilon^k}{N}\big{[}\frac{1}{\gamma}\langle q^{k+1/2}-p^k,p^*-p^k\rangle+\frac{1}{\gamma}\langle q^{k}-q^{k+1/2},p^*-p^k\rangle\big{]}\nonumber\\
&\leq&\frac{\epsilon^k}{N}\bigg{\{}\frac{1}{2\gamma}\big{[}\|p^*-p^k\|^2-\|p^*-q^{k+1/2}\|^2+\|q^{k+1/2}-p^k\|^2\big{]}\nonumber\\
&&+\frac{1}{\gamma}\|p^*-p^k\|\cdot\|q^{k}-q^{k+1/2}\|\bigg{\}}.\nonumber
\end{eqnarray}
Since $p^{k+1}=\mathcal{P}_\mu(q^{k+1/2})$ and property~\eqref{eq:Projecproperty4} of projection, we have $\|p^*-q^{k+1/2}\|\geq\|p^*-p^{k+1}\|$. Moreover, again using property~\eqref{eq:Projecproperty4} of projection, we obtain $\|q^{k}-q^{k+1/2}\|\leq\gamma\tau\|u^k-u^{k+1}\|$. Then~\eqref{eq:AL_dual_1} yields that
\begin{eqnarray}\label{eq:AL_dual_1-1}
\\
\frac{\epsilon^k}{N}\big{[}L_\gamma(u^{k},p^{*})-L_\gamma(u^k,p^{k})\big{]}&\leq&\frac{\epsilon^k}{N}\bigg{\{}\frac{1}{2\gamma}\big{[}\|p^*-p^k\|^2-\|p^*-p^{k+1}\|^2+\|q^{k+1/2}-p^k\|^2\big{]}\nonumber\\
&&+\tau\|p^*-p^k\|\cdot\|u^k-u^{k+1}\|\bigg{\}}.\nonumber
\end{eqnarray}
Since $\epsilon^{k+1}\leq\epsilon^k$,~\eqref{eq:AL_dual_1-1} yields that
\begin{eqnarray}\label{eq:AL_dual_1-2}
&&\frac{\epsilon^k}{N}\big{[}L_\gamma(u^{k},p^{*})-L_\gamma(u^k,p^{k})\big{]}\\
&\leq&\frac{1}{2\gamma N}\big{[}\epsilon^k\|p^*-p^k\|^2-\epsilon^{k+1}\|p^*-p^{k+1}\|^2\big{]}\nonumber\\
&&+\frac{\epsilon^k}{2\gamma N}\|q^{k+1/2}-p^k\|^2+\frac{\epsilon^k\tau}{N}\|p^*-p^k\|\cdot\|u^k-u^{k+1}\|\nonumber\\
&\leq&\frac{1}{2\gamma N}\big{[}\epsilon^k\|p^*-p^k\|^2-\epsilon^{k+1}\|p^*-p^{k+1}\|^2\big{]}\nonumber\\
&&+\frac{\epsilon^k}{\gamma N}[\|q^k-p^k\|^2+\|q^k-q^{k+1/2}\|^2]+\frac{\epsilon^k\tau}{N}\|p^*-p^k\|\cdot\|u^k-u^{k+1}\|\nonumber\\
&\leq&\frac{1}{2\gamma N}\big{[}\epsilon^k\|p^*-p^k\|^2-\epsilon^{k+1}\|p^*-p^{k+1}\|^2\big{]}\nonumber\\
&&+\frac{\epsilon^k}{\gamma N}\|q^k-p^k\|^2+\frac{\epsilon^0\gamma\tau^2}{ N}\|u^{k}-u^{k+1}\|^2+\frac{\epsilon^k\tau}{N}\|p^*-p^k\|\cdot\|u^k-u^{k+1}\|.\nonumber
\end{eqnarray}
Take expectation of $i(k)$ on both side of~\eqref{eq:AL_dual_1-2}, statement (ii) comes.\\
(iii) Note that $L_\gamma(u^{*},p^*)\geq L_\gamma(u^{*},p^k)$, we have that
\begin{eqnarray}\label{eq:AL_dual_3}
\\
\frac{\epsilon^k}{N}\big{[}L_\gamma(u^k,p^*)-L_\gamma(u^{*},p^*)\big{]}
&\leq&\frac{\epsilon^k}{N}\big{[}L_\gamma(u^k,p^*)-L_\gamma(u^{*},p^k)\big{]}\nonumber\\
&=&\frac{\epsilon^k}{N}\big{[}L_\gamma(u^{k},p^{*})-L_\gamma(u^k,p^{k})+L_\gamma(u^{k},p^{k})-L_\gamma(u^*,p^{k})\big{]}.\nonumber
\end{eqnarray}
Together statement (i) and (ii) of this lemma, we have
\begin{eqnarray}\label{eq:AL_dual_4}
\\
\frac{\epsilon^k}{N}\big{[}L_\gamma(u^k,p^*)-L_\gamma(u^{*},p^*)\big{]}&\leq&\Lambda^k(u^{k},p^k)-\mathbb{E}_{i(k)}\Lambda^{k+1}(u^{k+1},p^{k+1})\nonumber\\
&&+\epsilon^k\tau\big{(}\|q^k-p^*\|+\frac{1}{N}\|p^*-p^k\|\big{)}\cdot\mathbb{E}_{i(k)}\|u^{k}-u^{k+1}\|\nonumber\\
&&+\frac{\epsilon^k}{\gamma N}\|q^k-p^k\|^2+\frac{\epsilon^0\gamma\tau^2}{N}\mathbb{E}_{i(k)}\|u^{k}-u^{k+1}\|^2.\nonumber
\end{eqnarray}
To estimate the term $\mathbb{E}_{i(k)}\|u^{k}-u^{k+1}\|$, we consider~\eqref{eq:primal_VI} with  $u=u^k$, it follows that
\begin{eqnarray}\label{eq:AL_dual_5}
\\
&&\frac{1}{\epsilon^k}\langle\nabla K(u^{k+1})-\nabla K(u^k), u^{k+1}-u^k\rangle\nonumber\\
&\leq&\langle\nabla G(u^{k}),u^k-u^{k+1}\rangle+J(u^k)-J(u^{k+1})+\langle q^k,\Theta(u^k)-\Theta(u^{k+1})\rangle\nonumber\\
&\leq&\|\nabla G(u^{k})\|\cdot\|u^k-u^{k+1}\|+\|r^k\|\cdot\|u^k-u^{k+1}\|+\|q^k\|\cdot\|\Theta(u^k)-\Theta(u^{k+1})\|\nonumber\\
&&\qquad\qquad\qquad\qquad\qquad\qquad\qquad\qquad\qquad\qquad\qquad\qquad\qquad\mbox{(for $r^k\in\partial J(u^k)$)}\nonumber\\
&\leq&\big{(}\|\nabla G(u^{k})\|+\|r^k\|+\tau\|q^k\|\big{)}\|u^k-u^{k+1}\|.\nonumber
\end{eqnarray}
Since $\frac{1}{\epsilon^k}\langle\nabla K(u^{k+1})-\nabla K(u^k), u^{k+1}-u^k\rangle\geq\frac{\beta}{2\epsilon^k}\|u^{k+1}-u^k\|^2$, then we obtain
\begin{eqnarray}\label{eq:deltau}
\|u^{k+1}-u^k\|\leq\frac{2\epsilon^k}{\beta}\left(\|\nabla G(u^{k})\|+\|r^k\|+\tau\|q^k\|\right).
\end{eqnarray}
Therefore, inequality~\eqref{eq:AL_dual_4} yields that
\begin{eqnarray}
\\
\frac{\epsilon^k}{N}\big{[}L_\gamma(u^k,p^*)-L_\gamma(u^{*},p^*)\big{]}&\leq&\Lambda^k(u^{k},p^k)-\mathbb{E}_{i(k)}\Lambda^{k+1}(u^{k+1},p^{k+1})+\frac{\epsilon^k}{\gamma N}\|q^k-p^k\|^2\nonumber\\
&&+\frac{\epsilon^0\gamma\tau^2}{N}\mathbb{E}_{i(k)}\|u^{k}-u^{k+1}\|^2+h_1(u^k,p^k)(\epsilon^k)^2.\nonumber
\end{eqnarray}
Here comes the results.
\end{proof}
Next convergence theorem is derived by Proposition~\ref{proposition2} and Robbins-Siegmund Inequality~\cite{RS}.
\begin{theorem}\label{thm:convergence}{\bf(Almost surely convergence)}
\noindent Suppose Assumptions~\ref{assump1} and~\ref{assump2} hold, moreover if the nonincreasing sequence $\{\epsilon^k\}$ satisfies condition~\eqref{para-choice-convergence}, then
\begin{itemize}
\item[{\rm(i)}] $\sum\limits_{k=0}\limits^{+\infty}\epsilon^k\big{(}L_{\gamma}(u^k,p^*)-L_{\gamma}(u^*,p^*)\big{)}<+\infty\qquad\mbox{a.s.}$;
\item[{\rm(ii)}] Every cluster point of $\{u^{k}\}$ is almost surely an optimal solution of (P).
\end{itemize}
\end{theorem}
\begin{proof}
(i) Recalling statement (iii) of Lemma~\ref{lemma:bound2}, we have
\begin{eqnarray}
\\
\frac{\epsilon^k}{N}\big{[}L_\gamma(u^k,p^*)-L_\gamma(u^{*},p^*)\big{]}
&\leq&\Lambda^k(u^{k},p^k)-\mathbb{E}_{i(k)}\Lambda^{k+1}(u^{k+1},p^{k+1})+\frac{\epsilon^k}{\gamma N}\|q^k-p^k\|^2\nonumber\\
&&+\frac{\epsilon^0\gamma\tau^2}{N}\mathbb{E}_{i(k)}\|u^{k}-u^{k+1}\|^2+h_1(u^k,p^k)(\epsilon^k)^2.\nonumber
\end{eqnarray}
Then
\begin{eqnarray}
\\
\mathbb{E}_{i(k)}\Lambda^{k+1}(u^{k+1},p^{k+1})
&\leq&\Lambda^k(u^{k},p^k)+\frac{\epsilon^k}{\gamma N}\|q^k-p^k\|^2+\frac{\epsilon^0\gamma\tau^2}{N}\mathbb{E}_{i(k)}\|u^{k}-u^{k+1}\|^2\nonumber\\
&&+h_1(u^k,p^k)(\epsilon^k)^2-\frac{\epsilon^k}{N}\big{[}L_\gamma(u^k,p^*)-L_\gamma(u^{*},p^*)\big{]}.\nonumber
\end{eqnarray}
Since $\{u^k\}$ is almost surely bounded and $\{p^k\}$ is bounded, from Assumption~\ref{assump1}, we have $q^k$, $r^k$ and $\nabla G(u^k)$ are almost surely bounded. Then $h_1(u^k,p^k)$ is almost surely bounded. Therefore, $\sum\limits_{k=0}\limits^{+\infty}h_1(u^k,p^k)(\epsilon^k)^2<+\infty$, $\mbox{a.s.}$. From (i) of Propostion~\ref{proposition2}, we have $\sum\limits_{k=0}\limits^{+\infty}\mathbb{E}_{i(k)}\|u^k-u^{k+1}\|^2<+\infty$ $\mbox{a.s.}$ and $\sum\limits_{k=0}\limits^{+\infty}\epsilon^k\|q^k-p^{k}\|^2<+\infty$ $\mbox{a.s.}$. Noted $\frac{\epsilon^k}{N}\big{(}L_\gamma(u^k,p^*)-L_\gamma(u^{*},p^*)\big{)}\geq0$, again using Robbins-Siegmund Lemma~\cite{RS}, we get
\begin{eqnarray}\label{eq:covergence}
\sum\limits_{k=0}\limits^{+\infty}\epsilon^k\big{(}L_{\gamma}(u^k,p^*)-L_{\gamma}(u^*,p^*)\big{)}<+\infty\qquad\mbox{a.s.}.
\end{eqnarray}
(ii) Given that $L_\gamma(\cdot,p^*)$ is Lipschitz on every bounded set, combining~\eqref{eq:deltau} and (ii) of Proposition~\ref{proposition2}, we conclude from Lemma 4 of~\cite{CohenZ} that $\lim\limits_{k\rightarrow+\infty}L_{\gamma}(u^k,p^*)=L_\gamma(u^*,p^*)$ $\mbox{a.s.}$.\\
Let $\Omega_0$ denote the subset such that $\{u^k\}$ is not bounded, and let $\Omega_1$ denote the subset for which~\eqref{eq:covergence} does not hold: $\mathbb{P}(\Omega_0\cup\Omega_1)=0$. Pick some $\omega\notin\Omega_0\cup\Omega_1$. Since the sequence $\{u^k\}$ is bounded, it has cluster point and for each such cluster point $\bar{u}(\omega)$, it holds $\bar{u}(\omega)\in\arg\min\limits_{u\in\mathbf{U}} L_\gamma(u,p^*)$ from l.s.c. of $L_\gamma(\cdot,p^*)$. The stability of $L_\gamma$ implies that every $\bar{u}(\omega)$ is a solution of (P).  (See~\cite{CohenZ})
\end{proof}

\section{Almost surely convergence rate analysis}\label{ARate}
\begin{definition}
Given the desired accuracy $\varepsilon>0$, the primal point $u_\varepsilon\in\mathbf{U}$ is an $\varepsilon$-optimal solution for (P) if it satisfies\\
\indent Bound on primal suboptimality:
\begin{eqnarray}
&&\qquad|(G+J)(u_\varepsilon)-(G+J)(u^*)|\leq\varepsilon;\label{eq:suboptimality}
\end{eqnarray}
\indent Bound on feasibility:
\begin{eqnarray}
&&\qquad\qquad dist_{-\mathbf{C}}\big{(}\Theta(u_\varepsilon)\big{)}\leq\varepsilon.\label{eq:infeasibility}
\end{eqnarray}
\end{definition}
This section devotes to analyse the convergence rate for SPDC of (P). Before the convergence rate analysis, we again study the augmented Lagrangian function. Define the value function associated with (P$_1$) as
\begin{equation}\label{eq:valuefunction}
     v(\zeta)=\min\{(G+J)(u):\Theta(u)-\xi=\zeta,u\in\mathbf{U},\xi\in-\mathbf{C}\}, \forall\zeta\in\RR^m.
\end{equation}
From the convexity of $G+J$ and the $\mathbf{C}$-convexity of $\Theta$, it easy to show $v(\zeta)$ is convex in $\zeta$.
Then we write following problem
\begin{equation*}
\begin{array}{lll}
      \mbox{(P$_2$)}        &\min      & v(\zeta) \\
                            &{s.t}     & \zeta=0,
\end{array}
\end{equation*}
which is equivalent to (P) and (P$_1$).
The augmented Lagrangian function of (P$_2$) is
\begin{equation}\label{L_zeta}
\ell_\gamma(\zeta,p)=v(\zeta)+\langle p,\zeta\rangle+\frac{\gamma}{2}\|\zeta\|^2.
\end{equation}
The dual function for (P$_2$) is $\overline{\psi}_\gamma(p)=\min_{\zeta\in\RR^m}\ell_\gamma(\zeta,p)$. We also show that $\overline{\psi}_\gamma(p)$ is coincide with the dual function $\psi_\gamma(p)$ for (P).
Therefore, the optimal multiplier $p^*\in\mathbf{C}^*$ achieves the maximum of $\bar{\psi}_\gamma(p)$. From the stability of augmented Lagrangian of (P$_2$), we have $\hat{\zeta}_{p^*}=\arg\min \ell_\gamma(\zeta,p^*)=0$ and $\ell_\gamma(\hat{\zeta}_{p^*},p^*)=L_\gamma(u^*,p^*)=\psi_\gamma(p^*)$. From~\eqref{L_zeta}, $\ell_\gamma(\zeta,p)$ is strongly convex with constant $\gamma$, it yields
\begin{equation}\label{L_zeta_stronglyconvex}
\ell_\gamma(\zeta,p^*)-\ell_\gamma(\hat{\zeta}_{p^*},p^*)\geq\langle\nabla\ell_\gamma(\hat{\zeta}_{p^*},p^*),\zeta-\hat{\zeta}_{p^*}\rangle+\frac{\gamma}{2}\|\zeta-\hat{\zeta}_{p^*}\|^2\geq\frac{\gamma}{2}\|\zeta-\hat{\zeta}_{p^*}\|^2=\frac{\gamma}{2}\|\zeta\|^2.
\end{equation}
The following lemma furnishes the important character of augmented Lagrangian which will be used to derive the convergence rate.
\begin{lemma}\label{lemma:rate} {\bf(Inequalities for primal suboptimality and feasibility)}
\noindent Suppose Assumption~\ref{assump1} holds, then for any given $u\in\mathbf{U}$ and $p^*\in\mathbf{P}^*$, we have
\begin{eqnarray*}
\mbox{{\rm(i)}}\quad&&\|\Theta(u)-\xi(u,p^*)\|^2\leq\frac{2}{\gamma}[L_\gamma(u,p^*)-L_\gamma(u^*,p^*)]\quad\mbox{with}\quad\xi(u,p)=\arg\min_{\xi\in-\mathbf{C}}\overline{L}_{\gamma}(u,\xi,p);\\
\mbox{{\rm(ii)}}\quad&&|(G+J)(u)-(G+J)(u^*)|\leq[L_\gamma(u,p^*)-L_\gamma(u^*,p^*)]+\mu_0\sqrt{\frac{2}{\gamma}[L_\gamma(u,p^*)-L_\gamma(u^*,p^*)]};\\
\mbox{{\rm(iii)}}\quad&&\|\Pi\big{(}\Theta(u)\big{)}\|^2\leq\frac{2}{\gamma}[L_\gamma(u,p^*)-L_\gamma(u^*,p^*)].
\end{eqnarray*}
\end{lemma}
\begin{proof}
(i) Letting $\xi(u,p)=\arg\min\limits_{\xi\in-\mathbf{C}}\overline{L}_\gamma(u,\xi,p)$ and $\zeta(u,p)=\Theta(u)-\xi(u,p)$, $\forall u\in\mathbf{U}, p\in\mathbf{C}^*\cap \mathfrak{B}_\mu$. Observing that $(G+J)(u)\geq v(\zeta(u,p))$ due to the definition~\eqref{eq:valuefunction}, then
\begin{eqnarray}
L_\gamma(u,p)&=&\min_{\xi\in-\mathbf{C}}\overline{L}_\gamma(u,\xi,p)\\
             &=&(G+J)(u)+\langle p,\Theta(u)-\xi(u,p)\rangle+\frac{\gamma}{2}\|\Theta(u)-\xi(u,p)\|^2\nonumber\\
             &=&(G+J)(u)+\langle p,\zeta(u,p)\rangle+\frac{\gamma}{2}\|\zeta(u,p)\|^2\nonumber\\
             &\geq& v(\zeta(u,p))+\langle p,\zeta(u,p)\rangle+\frac{\gamma}{2}\|\zeta(u,p)\|^2\nonumber\\
             &\geq& \ell_\gamma(\zeta(u,p),p).\nonumber
\end{eqnarray}
This inequality implies
\begin{eqnarray}
L_\gamma(u,p^*)-L_\gamma(u^*,p^*)&\geq&\ell_\gamma(\zeta(u,p^*),p^*)-\ell_\gamma(\hat{\zeta}_{p^*},p^*)\\
                                 &\geq& \frac{\gamma}{2}\|\zeta(u,p^*)\|^2\qquad\mbox{(by~\eqref{L_zeta_stronglyconvex})}\nonumber\\
                                 &=& \frac{\gamma}{2}\|\Theta(u)-\xi(u,p^*)\|^2.\nonumber
\end{eqnarray}
(ii) From the definition of $L_\gamma(u,p)$ and $\xi(u,p)$, we have that
\begin{eqnarray*}
L_\gamma(u,p^*)-L_\gamma(u^*,p^*)&=&(G+J)(u)-(G+J)(u^*)+\langle p^*,\Theta(u)-\xi(u,p^*)\rangle+\frac{\gamma}{2}\|\Theta(u)-\xi(u,p^*)\|^2\nonumber\\
                                           &\geq&(G+J)(u)-(G+J)(u^*)+\langle p^*,\Theta(u)-\xi(u,p^*)\rangle\nonumber\\
                                           &\geq&(G+J)(u)-(G+J)(u^*)-\|p^*\|\cdot\|\Theta(u)-\xi(u,p^*)\|\nonumber\\
                                           &\geq&(G+J)(u)-(G+J)(u^*)-\mu_0\|\Theta(u)-\xi(u,p^*)\|.\qquad\mbox{(since $\|p^*\|\leq\mu_0$)}
\end{eqnarray*}
From statement (i) of this lemma, it follows that
\begin{equation}\label{rate1}
(G+J)(u)-(G+J)(u^*)\leq[L_\gamma(u,p^*)-L_\gamma(u^*,p^*)]+\mu_0\sqrt{\frac{2}{\gamma}[L_\gamma(u,p^*)-L_\gamma(u^*,p^*)]}.
\end{equation}
From the righthand side of saddle point inequality of $L_\gamma$, we obtain that
\begin{eqnarray*}
(G+J)(u^*)\leq(G+J)(u)+\langle p^*,\Theta(u)-\xi(u,p^*)\rangle+\frac{\gamma}{2}\|\Theta(u)-\xi(u,p^*)\|^2.
\end{eqnarray*}
Again use statement (i), consequently,
\begin{eqnarray*}
(G+J)(u)-(G+J)(u^*)&\geq&-\|p^*\|\cdot\|\Theta(u)-\xi(u,p^*)\|-\frac{\gamma}{2}\|\Theta(u)-\xi(u,p^*)\|^2\nonumber\\
                             &\geq&-[L_\gamma(u,p^*)-L_\gamma(u^*,p^*)]-\|p^*\|\sqrt{\frac{2}{\gamma}[L_\gamma(u,p^*)-L_\gamma(u^*,p^*)]}\nonumber\\
                             &\geq&-[L_\gamma(u,p^*)-L_\gamma(u^*,p^*)]-\mu_0\sqrt{\frac{2}{\gamma}[L_\gamma(u,p^*)-L_\gamma(u^*,p^*)]}.
\end{eqnarray*}
(iii) Since $\xi(u,p^*)\in-\mathbf{C}$, then
\begin{eqnarray}
\|\Theta(u)-\xi(u,p^*)\|^2&\geq& \|\Theta(u)-\Pi_{-\mathbf{C}}\big{(}\Theta(u)\big{)}\|^2\\
                          &=&\|\Pi\big{(}\Theta(u)\big{)}+\Pi_{-\mathbf{C}}\big{(}\Theta(u)\big{)}-\Pi_{-\mathbf{C}}\big{(}\Theta(u)\big{)}\|^2\nonumber\\
                          &&\qquad\qquad\mbox{(by property of projection~\eqref{eq:Projecproperty5})}\nonumber\\
                          &=&\|\Pi\big{(}\Theta(u)\big{)}\|^2.\nonumber
\end{eqnarray}
Together with statement (i) of this proposition, we have $\|\Pi\big{(}\Theta(u)\big{)}\|^2\leq\frac{2}{\gamma}[L_\gamma(u,p^*)-L_\gamma(u^*,p^*)]$.
\end{proof}
Noted that the augmented Lagrangian has the stability or it furnishes an exact penalty function for (P), i.e.
the solution of the following exact penalty problem is also one solution for (P).
\begin{equation}\label{eq:stability}
\begin{array}{lll}
      \mbox{(P$_3$)}&\min_{u\in\mathbf{U}} &L_\gamma(u,p^*).
\end{array}
\end{equation}
\indent Now we start the convergence rate analysis with above exact penalty problem (P$_3$) and the primal suboptimality and feasibility for problem (P).\\
\begin{theorem}\label{thm:rate}{\bf(Almost surely convergence rate)}
\noindent Let $\{u^k\}$ be the sequence generated from SPDC. Assume $\epsilon^k=\frac{1}{k^{\alpha}}$, $\frac{1}{2}<\alpha<1$, and $\bar{u}_t=\frac{\sum_{k=0}^{t}\epsilon^k u^k}{\sum_{k=0}^{t}\epsilon^k}$. Then under Assumption~\ref{assump1}, we have $d_1>0$ and
\begin{eqnarray*}
\mbox{{\rm(i)}}\quad&&L_\gamma(\bar{u}_{t},p^*)-L_\gamma(u^*,p^*)\leq\frac{d_1}{t^{1-\alpha}},\;\mbox{a.s.}\\
\mbox{{\rm(ii)}}\quad&&|(G+J)(\bar{u}_{t})-(G+J)(u^*)|\leq\frac{d_1}{t^{1-\alpha}}+\mu_0\sqrt{\frac{2d_1}{\gamma t^{1-\alpha}}},\;\mbox{a.s.}\\
\mbox{{\rm(iii)}}\quad&&dist_{-\mathbf{C}}\big{(}\Theta(\bar{u}_{t})\big{)}\leq\sqrt{\frac{2d_1}{\gamma t^{1-\alpha}}},\;\mbox{a.s.}.
\end{eqnarray*}
\end{theorem}
\begin{proof}
(i) From statement (i) of Theorem~\ref{thm:convergence} we have $d_1>0$ and
\begin{equation}
\sum\limits_{k=0}\limits^{t}\epsilon^k\big{(}L_{\gamma}(u^k,p^*)-L_{\gamma}(u^*,p^*)\big{)}<\sum\limits_{k=0}\limits^{+\infty}\epsilon^k\big{(}L_{\gamma}(u^k,p^*)-L_{\gamma}(u^*,p^*)\big{)}\leq d_1<+\infty\qquad\mbox{a.s.}
\end{equation}
By the definition $\bar{u}_t$, the convexity of $L_\gamma(u,p)$ for $u$ and $\epsilon^k=\frac{1}{k^{\alpha}}$, $\frac{1}{2}<\alpha\leq 1$, it follows that
\begin{equation}
L_\gamma(\bar{u}_{t},p^*)-L_\gamma(u^*,p^*)\leq\frac{\sum\limits_{k=0}\limits^{t}\epsilon^k\big{(}L_{\gamma}(u^k,p^*)-L_{\gamma}(u^*,p^*)}{\sum\limits_{k=0}\limits^{t}\epsilon^k}\leq\frac{d_1}{\sum\limits_{k=0}\limits^{t}\epsilon^k}\leq\frac{d_1}{t^{1-\alpha}}\qquad\mbox{a.s.}
\end{equation}
(ii) Together statement (ii) of Lemma~\ref{lemma:rate} and statement (i) of this theorem, we have that
\begin{equation}
|(G+J)(\bar{u}_{t})-(G+J)(u^*)|\leq\frac{d_1}{t^{1-\alpha}}+\mu_0\sqrt{\frac{2d_1}{\gamma t^{1-\alpha}}}\qquad\mbox{a.s.}
\end{equation}
(iii) From the statement (iii) of Lemma~\ref{lemma:rate} and statement (i) of this theorem, it follows
\begin{equation}
dist_{-\mathbf{C}}\big{(}\Theta(\bar{u}_{t})\big{)}=\|\Theta(\bar{u}_{t})-\Pi_{-\mathbf{C}}\big{(}\Theta(\bar{u}_{t})\big{)}\|=\|\Pi\big{(}\Theta(\bar{u}_{t})\big{)}\|\leq\sqrt{\frac{2d_1}{\gamma t^{1-\alpha}}}\qquad\mbox{a.s.}
\end{equation}
\end{proof}
\begin{remark}
Theorem~\ref{thm:rate} shows $\bar{u}_t$ produced by Algorithm SPDC almost surely satisfies the suboptimality and feasibility condition~\eqref{eq:suboptimality} and~\eqref{eq:infeasibility}. The convergence rate is $O(1/t^{\frac{1-\alpha}{2}})$ in this case, i.e. given a desired accuracy $\varepsilon>0$ the SPDC produces an $\varepsilon$-solution of (P) in $O(1/\varepsilon^{\frac{2}{1-\alpha}})$ iterations with probability $1$.
\end{remark}
\indent Next section we will propose an expected convergence rate generated by algorithm SPDC.
\section{Expected Convergence rate}\label{ERate}
In this section we prove expected convergence rate for SPDC. Before proceeding, we need following lemma for the bifunction value at $(u^k,q^k)$.
\begin{lemma}\label{lemma:bound3}{\bf(Monotonicity of expected distance like function)}
Suppose Assumption~\ref{assump1} and~\ref{assump2} hold, $\{(u^k,p^k)\}$ is generated by SPDC, the parameter sequence $\{\epsilon^k\}$ satisfies condition~\eqref{para-choice-convergence}. Then for any $t>k>0$ and $p\in\mathbf{C}^*\cap\mathfrak{B}_\mu$, it holds that
\begin{eqnarray*}
\frac{\epsilon^k}{N}\mathbb{E}_{\mathcal{F}_{t}}\big{[}L(u^k,p)-L(u^*,q^{k})\big{]}\leq\mathbb{E}_{\mathcal{F}_{t}}\bigg{[}D(u^*,u^k)-D(u^*,u^{k+1})+\frac{\epsilon^k}{2\gamma N}\|p-p^k\|^2-\frac{\epsilon^{k+1}}{2\gamma N}\|p-p^{k+1}\|^2+h_2(u^k,q^k)(\epsilon^k)^2\bigg{]};
\end{eqnarray*}
where\\
$h_2(u^k,q^k)=\frac{2}{\beta}\big{(}\|\nabla G(u^{k})\|+\|r^k\|+\frac{(N+1)\tau}{N}\|q^k\|+\frac{\tau}{N}\mu\big{)}(\|\nabla G(u^{k})\|+\|r^k\|+\tau\|q^k\|)$ and $p$ could possibly be random.\\
\end{lemma}
\begin{proof}
From (i) of Lemma~\ref{lemma:bound1}, we obtain
\begin{eqnarray}\label{eq:bound2_2_1}
&&\frac{\epsilon^k}{N}\big{[}L(u^{k}, q^{k})-L(u^*, q^{k})\big{]}\\
&\leq&\big{[}D(u^*,u^k)+\epsilon^k\big{(}L(u^{k},p^*)-L(u^*,p^*)\big{)}\big{]}\nonumber\\
&&-\mathbb{E}_{i(k)}\big{[}D(u^*,u^{k+1})+\epsilon^k\big{(}L(u^{k+1},p^*)-L(u^*,p^*)\big{)}\big{]}\nonumber\\
&&+\epsilon^k\mathbb{E}_{i(k)}\langle q^k-p^*,\Theta(u^{k})-\Theta(u^{k+1})\rangle-\frac{\beta-\epsilon^kB_G}{2}\mathbb{E}_{i(k)}\|u^k-u^{k+1}\|^2\nonumber\\
&=&\big{[}D(u^*,u^k)-\mathbb{E}_{i(k)}D(u^*,u^{k+1})\big{]}+\epsilon^k\mathbb{E}_{i(k)}\big{[}(G+J)(u^{k})-(G+J)(u^{k+1})\big{]}\nonumber\\
&&+\epsilon^k\mathbb{E}_{i(k)}\langle q^k,\Theta(u^{k})-\Theta(u^{k+1})\rangle-\frac{\beta-\epsilon^kB_G}{2}\mathbb{E}_{i(k)}\|u^k-u^{k+1}\|^2\nonumber\\
&\leq&\big{[}D(u^*,u^k)-\mathbb{E}_{i(k)}D(u^*,u^{k+1})\big{]}+\epsilon^k\big{(}\|\nabla G(u^{k})\|+\|r^k\|+\tau\|q^k\|\big{)}\mathbb{E}_{i(k)}\|u^k-u^{k+1}\|\nonumber\\
&&-\frac{\beta-\epsilon^kB_G}{2}\mathbb{E}_{i(k)}\|u^k-u^{k+1}\|^2.\nonumber
\end{eqnarray}
Taking expectation with respect to $\mathcal{F}_{t}$, $t>k$, for inequality~\eqref{eq:bound2_2_1}, we obtain that
\begin{eqnarray}\label{eq:bound2_2_2}
&&\frac{\epsilon^k}{N}\mathbb{E}_{\mathcal{F}_{t}}\big{[}L(u^{k}, q^{k})-L(u^*, q^{k})\big{]}\\
&\leq&\mathbb{E}_{\mathcal{F}_{t}}\bigg{[}D(u^*,u^k)-D(u^*,u^{k+1})+\epsilon^k\big{(}\|\nabla G(u^{k})\|+\|r^k\|+\tau\|q^k\|\big{)}\|u^k-u^{k+1}\|\nonumber\\
&&-\frac{\beta-\epsilon^kB_G}{2}\|u^k-u^{k+1}\|^2\bigg{]}.\nonumber
\end{eqnarray}
Taking expectation with respect to $\mathcal{F}_{t}$ for (ii) of Lemma~\ref{lemma:bound1}, it follows that, $\forall p\in\mathbf{C}^*\cap\mathfrak{B}_{\mu}$,
\begin{eqnarray}\label{eq:bound2_2_3}
&&\frac{\epsilon^k}{N}\mathbb{E}_{\mathcal{F}_{t}}\big{[}L(u^k,p)-L(u^{k},q^{k})\big{]}\\
&\leq&\mathbb{E}_{\mathcal{F}_{t}}\bigg{[}\frac{\epsilon^k}{2\gamma N}\|p-p^k\|^2-\frac{\epsilon^{k+1}}{2\gamma N}\|p-p^{k+1}\|^2\nonumber\\
&&-\frac{\epsilon^k}{N}\langle q^k-p,\Theta(u^{k})-\Theta(u^{k+1})\rangle+\frac{\gamma\tau^2\epsilon^k}{2N}\|u^{k}-u^{k+1}\|^2-\frac{\epsilon^k}{2\gamma N}\|q^k-p^{k}\|^2\bigg{]}\nonumber\\
&\leq&\mathbb{E}_{\mathcal{F}_{t}}\bigg{[}\frac{\epsilon^k}{2\gamma N}\|p-p^k\|^2-\frac{\epsilon^{k+1}}{2\gamma N}\|p-p^{k+1}\|^2\nonumber\\
&&+\frac{\epsilon^k\tau}{N}[\|q^k\|+\mu]\|u^{k}-u^{k+1}\|+\frac{\gamma\tau^2\epsilon^k}{2N}\|u^{k}-u^{k+1}\|^2\bigg{]}\nonumber
\end{eqnarray}
Summing~\eqref{eq:bound2_2_2} and~\eqref{eq:bound2_2_3}, it follows that
\begin{eqnarray*}
&&\frac{\epsilon^k}{N}\mathbb{E}_{\mathcal{F}_{t}}\big{[}L(u^k,p)-L(u^*,q^{k})\big{]}\\
&\leq&\mathbb{E}_{\mathcal{F}_{t}}\bigg{[}D(u^*,u^k)-D(u^*,u^{k+1})+\frac{\epsilon^k}{2\gamma N}\|p-p^k\|^2-\frac{\epsilon^{k+1}}{2\gamma N}\|p-p^{k+1}\|^2\nonumber\\
&&+\epsilon^k\big{[}\|\nabla G(u^{k})\|+\|r^k\|+\frac{(N+1)\tau}{N}\|q^k\|+\frac{\tau}{N}\mu\big{]}\|u^k-u^{k+1}\|\nonumber\\
&&-\frac{N\beta-\epsilon^k(NB_G+\gamma\tau^2)}{2N}\|u^k-u^{k+1}\|^2\bigg{]}\nonumber\\
&\leq&\mathbb{E}_{\mathcal{F}_{t}}\bigg{[}D(u^*,u^k)-D(u^*,u^{k+1})+\frac{\epsilon^k}{2\gamma N}\|p-p^k\|^2-\frac{\epsilon^{k+1}}{2\gamma N}\|p-p^{k+1}\|^2\nonumber\\
&&+\epsilon^k\big{[}\|\nabla G(u^{k})\|+\|r^k\|+\frac{(N+1)\tau}{N}\|q^k\|+\frac{\tau}{N}\mu\big{]}\|u^k-u^{k+1}\|\bigg{]}\qquad\qquad\qquad\qquad\qquad\qquad\quad\mbox{(since $\epsilon^k$ satisfy~\eqref{para-choice-convergence})}\nonumber\\
&\leq&\mathbb{E}_{\mathcal{F}_{t}}\bigg{[}D(u^*,u^k)-D(u^*,u^{k+1})+\frac{\epsilon^k}{2\gamma N}\|p-p^k\|^2-\frac{\epsilon^{k+1}}{2\gamma N}\|p-p^{k+1}\|^2+h_2(u^k,q^k)(\epsilon^k)^2\bigg{]}\qquad\mbox{(from~\eqref{eq:deltau})}
\end{eqnarray*}
Noted that the above derivation is valid even though $p$ is a random vector. Then, we have the claimed results.
\end{proof}
For the sequence $\{(u^k,p^k)\}$ generated from Algorithm SPDC, and any $t>0$ we define the average sequence $\bar{u}_{t}=\frac{\sum_{k=0}^{t}\epsilon^k u^k}{\sum_{k=0}^{t}\epsilon^k}$ and $\bar{p}_{t}=\frac{\sum_{k=0}^{t}\epsilon^k q^k}{\sum_{k=0}^{t}\epsilon^k}$. Then we have the following proposition.
\begin{proposition}\label{thm:rate_ergodic} {\bf(Global estimate of expected bifunction values in $(\bar{u}_t,\bar{p}_t)$)}
\noindent Let Assumption~\ref{assump1} and~\ref{assump2} hold, $\epsilon^k=\frac{1}{k^\alpha}$. Then we have $d_2>0$ such that
$$\mathbb{E}_{\mathcal{F}_{t}}\big{[}L(\bar{u}_t,p)-L(u^*,\bar{p}_t)\big{]}\leq\frac{Nd_2}{t^{1-\alpha}}, \forall p\in\mathbf{C}^*\cap\mathfrak{B}_\mu,$$
where $p$ could possibly be random.
\end{proposition}
\begin{proof}
For any given integer $k$, $u^k$ is almost surely bounded and $p^k$ is bounded, then $q^k$ is also almost surely bounded, there exists $\delta\geq\mathbb{E}_{\mathcal{F}_{t}}h_2(u^k,q^k)$. Then from Lemma~\eqref{lemma:bound3}, we obtain
\begin{eqnarray*}
\frac{\epsilon^k}{N}\mathbb{E}_{\mathcal{F}_{t}}\big{[}L(u^k,p)-L(u^*,q^{k})\big{]}\leq\mathbb{E}_{\mathcal{F}_{t}}\bigg{[}D(u^*,u^k)-D(u^*,u^{k+1})+\frac{\epsilon^k}{2\gamma N}\|p-p^k\|^2-\frac{\epsilon^{k+1}}{2\gamma N}\|p-p^{k+1}\|^2+\delta(\epsilon^k)^2\bigg{]}.
\end{eqnarray*}
Summing it over $k=0$ through $t$, we have $\forall p\in\mathbf{C}^*\cap\mathfrak{B}_{\mu}$,
\begin{eqnarray*}
\sum_{k=0}^{t}\epsilon^k\mathbb{E}_{\mathcal{F}_{t}}\big{[}L(u^k,p)-L(u^*,q^k)\big{]}&\leq& N\mathbb{E}_{\mathcal{F}_{t}}[D(u^*,u^0)+\frac{\epsilon^k}{2\gamma N}\|p-p^0\|^2+\delta\sum_{k=0}^{t}(\epsilon^k)^2]\\
&\leq&N\bigg{[}\frac{B}{2}\|u^0-u^*\|^2+\frac{2\epsilon^0\mu^2}{\gamma N}+\delta\sum_{k=0}^{t}(\epsilon^k)^2\bigg{]}.
\end{eqnarray*}
Moreover, since $\sum\limits_{k=0}\limits^{+\infty}(\epsilon^k)^2<+\infty$, then there exists $d_2>0$, such that $d_2\geq\frac{B}{2}\|u^0-u^*\|^2+\frac{2\epsilon^0\mu^2}{\gamma N}+\delta\sum_{k=0}^{t}(\epsilon^k)^2$, and we have
\begin{eqnarray}\label{rate2}
\sum_{k=0}^{t}\epsilon^k\mathbb{E}_{\mathcal{F}_{t}}\big{[}L(u^k,p)-L(u^*,q^k)\big{]}\leq Nd_2.
\end{eqnarray}
Another hand, from the definition of $\bar{u}_t$ and $\bar{p}_t$, we have $\bar{u}_t\in\mathbf{U}$ and $\bar{p}_t\in\mathbf{C}^*$. From the convexity of set $\mathbf{U}$, $\mathbf{C}^*$ and the function $L(u',p)-L(u,p')$ is convex in $u'$ and linear in $p'$, for all $p\in\mathbf{C}^*\cap\mathfrak{B}_\mu$, we have that
\begin{eqnarray}\label{eq:VI_rate1}
\mathbb{E}_{\mathcal{F}_{t}}\big{[}L(\bar{u}_t,p)-L(u^*,\bar{p}_t)\big{]}&\leq&\frac{1}{\sum_{k=0}^{t}\epsilon^k}\sum_{k=0}^{t}\epsilon^k\mathbb{E}_{\mathcal{F}_{t}}\big{[}L(u^k,p)-L(u^*,q^k)\big{]}\\
                                                 &\leq&\frac{Nd_2}{\sum_{k=0}^{t}\epsilon^k}.\nonumber\qquad\mbox{(by~\eqref{rate2})}
\end{eqnarray}
Finally, taking $\epsilon^k=\frac{1}{k^{\alpha}}$ with $\frac{1}{2}<\alpha< 1$, we conclude
\begin{eqnarray}\label{eq:VI_rate2}
\mathbb{E}_{\mathcal{F}_{t}}\big{[}L(\bar{u}_t,p)-L(u^*,\bar{p}_t)\big{]}\leq\frac{Nd_2}{\sum_{k=0}^{t}\frac{1}{k^\alpha}}\leq\frac{Nd_2}{t^{1-\alpha}}.
\end{eqnarray}
Noted that the above derivation is valid even though $p$ is a random vector.
\end{proof}
Next theorem provides the rate of expected feasibility and primal suboptimality.
\begin{theorem}\label{Co:2}{\bf(Convergence rate of expected feasibility and primal suboptimality)}
Suppose assumptions of Proposition~\ref{thm:rate_ergodic} hold, then the following assertion of convergence rate hold:
\begin{eqnarray*}
\mbox{{\rm(i)}}\quad&&\mathbb{E}_{\mathcal{F}_{t}}dist_{-\mathbf{C}}\big{(}\Theta(\bar{u}_{t})\big{)}\leq\frac{Nd_2}{\gamma t^{1-\alpha}};\\
\mbox{{\rm(ii)}}\quad&&-\frac{\mu_0 Nd_2}{t^{1-\alpha}}\leq \mathbb{E}_{\mathcal{F}_{t}}[(G+J)(\bar{u}_{t})-(G+J)(u^*)]\leq\frac{Nd_2}{t^{1-\alpha}}.
\end{eqnarray*}
\end{theorem}
\begin{proof}
(i) If $\mathbb{E}_{\mathcal{F}_{t}}\Pi\big{(}\Theta(\bar{u}_{t})\big{)}=0$, statement (i) is obviously. Otherwise, $\mathbb{E}_{\mathcal{F}_{t}}\Pi\big{(}\Theta(\bar{u}_{t})\big{)}\neq0$ i.e., there is set $\Omega_1$ such that $\mathbb{P}\big{\{}\omega\in\Omega_1\big{|}\Pi\big{(}\Theta(\bar{u}_{t}(\omega))\big{)}\neq 0\big{\}}>0$.\\
\indent Let $\hat{p}$ be a random vector:
\begin{eqnarray*}
\hat{p}(\omega)=\left\{
\begin{array}{cc}
0 & \omega\notin\Omega_1 \\
\frac{\mu\Pi\big{(}\Theta(\bar{u}_t(\omega))\big{)}}{\|\Pi\big{(}\Theta(\bar{u}_t(\omega))\big{)}\|} & \omega\in\Omega_1.
\end{array}
\right.
\end{eqnarray*}
Thus $\hat{p}\in\mathbf{C}^*\cap\mathfrak{B}_{\mu}$. Noted that $\hat{p}(\omega)=0$, $\omega\notin\Omega_1$. For all $\omega\in\Omega_1$, we have that
\begin{eqnarray}
\\
\langle \hat{p}(\omega),\Theta(\bar{u}_t(\omega))\rangle&=&\langle \frac{\mu\Pi\big{(}\Theta(\bar{u}_t(\omega))\big{)}}{\|\Pi\big{(}\Theta(\bar{u}_t(\omega))\big{)}\|},\Theta(\bar{u}_t(\omega))\rangle\nonumber\\
&=&\langle\frac{\mu\Pi\big{(}\Theta(\bar{u}_{t}(\omega))\big{)}}{\|\Pi\big{(}\Theta(\bar{u}_{t}(\omega))\big{)}\|},\Pi\big{(}\Theta(\bar{u}_{t}(\omega))\big{)}+\Pi_{-\mathbf{C}}\big{(}\Theta(\bar{u}_{t}(\omega))\big{)}\rangle\quad\mbox{(since~\eqref{eq:Projecproperty5})}\nonumber\\
&=&\mu\|\Pi\big{(}\Theta(\bar{u}_t(\omega))\big{)}\|.\qquad\qquad\qquad\qquad\qquad\qquad\qquad\quad\;\mbox{(since~\eqref{eq:Projecproperty6})}\nonumber
\end{eqnarray}
Then we have
\begin{equation}\label{pp}
\mathbb{E}_{\mathcal{F}_{t}}\langle\hat{p},\Theta(\bar{u}_t)\rangle=\mu\mathbb{E}_{\mathcal{F}_{t}}\|\Pi\big{(}\Theta(\bar{u}_t)\big{)}\|.
\end{equation}
Let $\{u^k\}$ be the sequence generated by SPDC. By Proposition~\ref{thm:rate_ergodic} with $p=\hat{p}$, we obtain that
\begin{eqnarray}\label{eq:inequal_rate1}
&&Nd_2/t^{1-\alpha}\\
&\geq&\mathbb{E}_{\mathcal{F}_{t}}\big{[}L(\bar{u}_t,\hat{p})-L(u^*,\bar{p}_t)\big{]}\nonumber\\
&=&\mathbb{E}_{\mathcal{F}_{t}}[(G+J)(\bar{u}_{t})-(G+J)(u^*)+\langle \hat{p},\Theta(\bar{u}_{t})\rangle-\langle\bar{p}_t,\Theta(u^*)\rangle]\nonumber\\
&\geq&\mathbb{E}_{\mathcal{F}_{t}}[(G+J)(\bar{u}_{t})-(G+J)(u^*)+\langle\hat{p},\Theta(\bar{u}_{t})\rangle]\nonumber\\
&&\qquad\qquad\qquad\qquad\qquad\qquad\qquad\;\mbox{(since $\bar{p}_t\in\mathbf{C}^*$ and $\Theta(u^*)\in-\mathbf{C}$.)}\nonumber\\
&=&\mathbb{E}_{\mathcal{F}_{t}}[(G+J)(\bar{u}_{t})-(G+J)(u^*)+\mu\|\Pi\big{(}\Theta(\bar{u}_{t})\big{)}\|]\nonumber
\end{eqnarray}
Take $u=\bar{u}_{t}$ in the right hand side of saddle point inequality~\eqref{saddle point:L}, we have
\begin{eqnarray}\label{eq:inequal_rate2}
(G+J)(\bar{u}_{t})-(G+J)(u^*)\geq-\|p^*\|\cdot\|\Pi\big{(}\Theta(\bar{u}_{t})\big{)}\|\geq-\mu_0\|\Pi\big{(}\Theta(\bar{u}_{t})\big{)}\|.
\end{eqnarray}
Together~\eqref{eq:inequal_rate1} and~\eqref{eq:inequal_rate2}, we have
\begin{equation}\label{eq:inequal_rate3}
\mathbb{E}_{\mathcal{F}_{t}}\|\Pi\big{(}\Theta(\bar{u}_{t})\big{)}\|\leq\frac{Nd_2}{ t^{1-\alpha}}\qquad\qquad\mbox{(since $\mu=\mu_0+1$)}
\end{equation}
Thus we obtain
\begin{equation}\label{eq:inequal_rate4}
\mathbb{E}_{\mathcal{F}_{t}}dist_{-\mathbf{C}}\big{(}\Theta(\bar{u}_{t})\big{)}=\mathbb{E}_{\mathcal{F}_{t}}\|\Theta(\bar{u}_{t})-\Pi_{-\mathbf{C}}\big{(}\Theta(\bar{u}_{t})\big{)}\|=\mathbb{E}_{\mathcal{F}_{t}}\|\Pi\big{(}\Theta(\bar{u}_{t})\big{)}\|\leq\frac{Nd_2}{ t^{1-\alpha}}
\end{equation}
Here comes statement (i).\\
(ii) Together~\eqref{eq:inequal_rate1}-\eqref{eq:inequal_rate3}, we have
\begin{equation}
-\frac{\mu_0Nd_2}{ t^{1-\alpha}}\leq\mathbb{E}_{\mathcal{F}_{t}}[(G+J)(\bar{u}_{t})-(G+J)(u^*)]\leq \frac{Nd_2}{ t^{1-\alpha}}.
\end{equation}
\end{proof}
\begin{remark}{\bf(High probability complexity bound for $\varepsilon$-solution)}
From Theorem~\ref{Co:2}, we immediately get
\begin{eqnarray*}
\mathbb{E}_{\mathcal{F}_{t}}\bigg{\{}|(G+J)(\bar{u}_{t})-(G+J)(u^*)|+dist_{-\mathbf{C}}\big{(}\Theta(\bar{u}_{t})\big{)}\bigg{\}}\leq\frac{(1+\mu)Nd_2}{ t^{1-\alpha}}
\end{eqnarray*}
from Markov's inequality~\cite{LiuWright1, LiuWright2, LuXiao}, it follows
\begin{eqnarray*}
&&\mathbb{P}\bigg{\{}|(G+J)(\bar{u}_{t})-(G+J)(u^*)|+dist_{-\mathbf{C}}\big{(}\Theta(\bar{u}_{t})\big{)}\geq\varepsilon\bigg{\}}\nonumber\\
&\leq&\varepsilon^{-1}\mathbb{E}_{\mathcal{F}_{t}}\bigg{\{}|(G+J)(\bar{u}_{t})-(G+J)(u^*)|+dist_{-\mathbf{C}}\big{(}\Theta(\bar{u}_{t})\big{)}\bigg{\}}\nonumber\\
&\leq&\frac{(1+\mu)Nd_2}{\varepsilon t^{1-\alpha}}\nonumber\\
&\leq&\eta,
\end{eqnarray*}
providing the following condition holds, i.e.,
\begin{eqnarray}
t\geq\left(\frac{(1+\mu)Nd_2}{\varepsilon\eta}\right)^{\frac{1}{1-\alpha}}.
\end{eqnarray}
\end{remark}
Theorem~\ref{thm:rate_ergodic} shows that SPDC has the expected convergence rate $O(1/t^{1-\alpha})$ in the
worst case. Given a desire accuracy $\varepsilon>0$ the SPDC method produces an $\varepsilon$-solution of (P) in $\left(\frac{(1+\mu)Nd_2}{\varepsilon\eta}\right)^{\frac{1}{1-\alpha}}$ iterations with probability $1-\eta$.
\section{Numerical Experiments}\label{sec:experiments}
\begin{figure}
\begin{minipage}[t]{0.45\linewidth}
\centering
\includegraphics[width=2.2in]{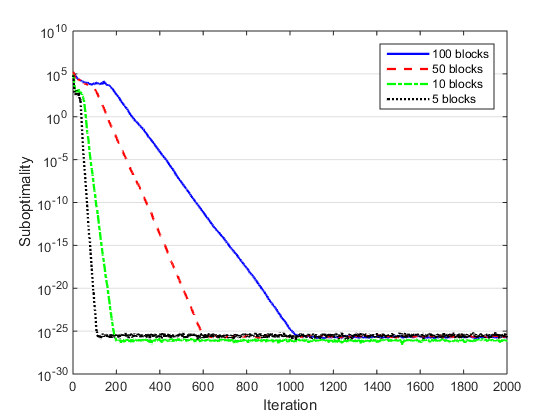}
\label{fig:side:1-b}
\end{minipage}
\begin{minipage}[t]{0.5\linewidth}
\centering
\includegraphics[width=2.2in]{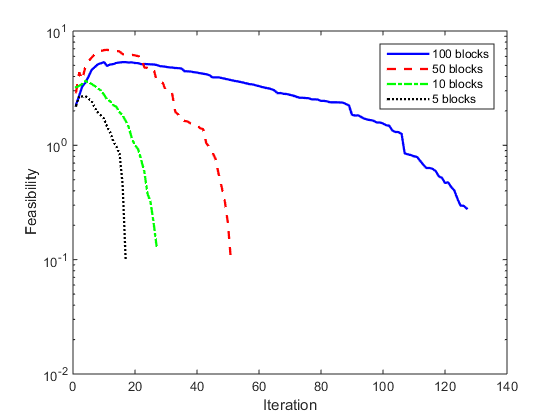}
\label{fig:side:1-c}
\end{minipage}
\caption{$m=1200$, $n=2000$, and $s=10$. The left-hand graph in each figure shows the number of blocks and plots suboptimality versus iteration count. The right-hand graph indicates the number of blocks and plots feasibility value versus iteration count.}
\label{fig:1}
\end{figure}

\begin{figure}
\begin{minipage}[t]{0.45\linewidth}
\centering
\includegraphics[width=2.1in]{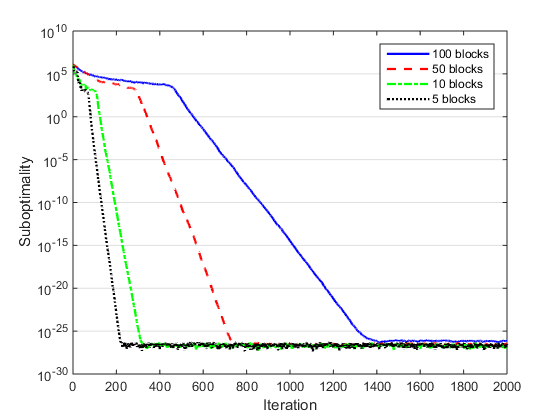}
\label{fig:side:2-b}
\end{minipage}
\begin{minipage}[t]{0.5\linewidth}
\centering
\includegraphics[width=2.1in]{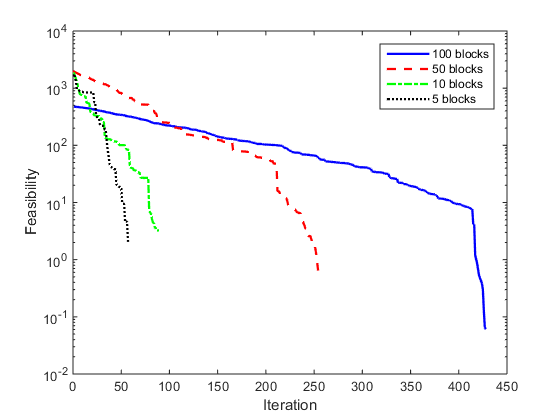}
\label{fig:side:2-c}
\end{minipage}
\caption{$m=3000$, $n=5000$, and $s=25$. The left-hand graph in each figure shows the number of blocks and plots suboptimality versus iteration count. The right-hand graph indicates the number of blocks and plots feasibility value versus iteration count.}
\label{fig:2}
\end{figure}

\indent In this section, we test the proposed SPDC method on solving the Elastic Net Support Vector Machine problem (EN-SVM)~\cite{ZouHastie}:
\begin{equation}\label{Prob:EN-SVM}
\begin{array}{lll}
\mbox{(EN-SVM):}  &\min       & \frac{1}{2}\|Au-b\|^2      \\
                  &\rm {s.t}  & \lambda\|u\|_1+(1-\lambda)\|u\|^2\leq\delta \\
                  &           & u\in\RR^n.
\end{array}
\end{equation}
\indent The elements of $A\in\RR^{m\times n}$ are selected i.i.d. from a Gaussian $\mathcal{N}(0,1)$ distribution. To construct a sparse true solution $u^*\in\RR^n$, given the dimension $n$ and sparsity $s$, we select $s$ entries of $u^*$ at random to be nonzero and $\mathcal{N}(0,1)$ normally distributed, and set the rest to zero. The measurement vector $b\in\RR^m$ is obtained by $b=Au^*$, $\lambda\in[0,1]$ and $\delta\in\RR_+$.\\
\indent We choose $\lambda=0.4$ and $\delta=\lambda\|u^*\|_1+(1-\lambda)\|u^*\|^2$ with $m=200$, $n=2000$, and $s=10$ in Figure~\ref{fig:1} and $m=500$, $n=5000$, and $s=25$ in Figure~\ref{fig:2}. It is obvious that the optimal value of both cases are zero. In both cases, we partition the variables into $5$, $10$, $50$ and $100$ blocks (i.e. $N=5,10,50,100$). Then in Figure~\ref{fig:1}, $n_i=400,200,40,20$. And in Figure~\ref{fig:2}, $n_i=1000,500,100,50$. In each iteration, we randomly choose one block $i$ of variables to update as follows:
\begin{eqnarray}
u_i^{k+1}&=&\arg\min\frac{1+(1-\lambda)2\epsilon^kq^k}{2\epsilon^k\lambda q^k}\bigg{\|}u_i-\frac{u_i^k+\epsilon^kA_i^{\top}(Au^k-b)}{1+2\epsilon^k(1-\lambda)q^k}\bigg{\|}^2+\|u_i\|_1;\\
p^{k+1}&=&\min\{\max[p^k+\gamma\lambda\|u^{k+1}\|_1+\gamma(1-\lambda)\|u^{k+1}\|^2,0],\mu\};
\end{eqnarray}
where $q^k=\max[p^k+\gamma\lambda\|u^{k}\|_1+\gamma(1-\lambda)\|u^{k}\|^2,0]$ and $\mu=\frac{\|b\|^2}{2\delta}+1$. Here the primal problem has a closed form solution as follows:
$$u_i^{k+1}=sign(r)\odot\max\bigg{\{}0,|r|-\mathbf{1}_{n_i}\cdot\frac{\epsilon^k\lambda q^k}{1+(1-\lambda)2\epsilon^kq^k}\bigg{\}},$$
where $r=\frac{u_i^k+\epsilon^kA_i^{\top}(Au^k-b)}{1+2\epsilon^k(1-\lambda)q^k}$ and $\odot$ denotes componentwise multiplication. Moreover, we select $\epsilon^k$ as $\epsilon^{k}=\frac{1}{1000+(k/1000)}$ in Figure~\ref{fig:1} and $\epsilon^{k}=\frac{1}{10000+(k/1000)}$ in Figure~\ref{fig:2}. We perform two experiments in MATLAB(R2011b) on a personal computer with an Intel Core i5-6200U CPUs (2.40GHz) and 8.00 GB of RAM.\\
\indent The left-hand graph in each figure shows the number of blocks and plots suboptimality versus iteration count. The right-hand graph indicates the number of blocks and plots feasibility value versus iteration count.


\section*{Acknowledgments}
We would like to acknowledge the assistance of volunteers in putting
together this example manuscript and supplement.

\bibliographystyle{siamplain}

\end{document}